\documentclass[12pt]{amsart}
\usepackage{header}
\usepackage{url}
\numberwithin{equation}{section}
\usepackage{cite}
\usepackage{subfiles} 
\begin{document}
\title{Directional Ballistic transport in Quantum Waveguides}
\author[Black]{Adam Black}
\address{Department of Mathematics\\
  University of California, Berkeley\\
 Berkeley, CA 94720}
\email[]{adamblack@berkeley.edu}
\thanks{This material is based upon work supported by the National Science Foundation under Award No. 2503339}
\author[Damanik]{David Damanik}
\address{Department of Mathematics\\
  Rice University\\
 Houston, TX, 77005}
\email[]{damanik@rice.edu }
\thanks{D.\ D.\ was supported in part by NSF grant DMS--2349919}
\author[Kuchment]{Peter Kuchment}
\address{Department of Mathematics\\
 Texas A\&M University\\
 College Station, TX 77843}
 \thanks{P.K. thanks the Arthur George and Mary Emolene Owen endowment for partial support}
\author[Malinovitch]{Tal Malinovitch}
\address{Department of Mathematics\\
  Rice University\\
 Houston, TX, 77005}
  \email[]{tal.malinovitch@rice.edu}
\author[Young]{Giorgio Young}
\address{Department of Mathematics\\
  The University of Wisconsin\\
 Madison, WI 53706}
  \email[]{gfyoung@wisc.edu}
  \thanks{G.Y.\ acknowledges the support of the National Science Foundation through grant DMS--2303363.}
\maketitle
\begin{abstract}
     We study the transport properties of Schr\"{o}dinger operators on $\bbR^d$ with potentials that are periodic in some directions and compactly supported in the others. Such systems are known to produce \emph{surface states} that are weakly confined near the support of the potential. We show that a natural set of surface states exhibits \emph{directional ballistic transport}, characterized by ballistic transport in the periodic directions and its absence in the others. To prove this, we develop a Floquet theory that captures the analytic variation of surface states. The main idea consists of reformulating the eigenvalue problem for surface states as a Fredholm problem via the Dirichlet-to-Neumann map. 
\end{abstract}

\section{Introduction}
Elliptic operators with coefficients that are periodic only in some coordinate directions appear naturally in many physical contexts, where they may describe \emph{waveguides}.
Such systems are interesting because they can support guided modes that propagate only in certain directions, unlike typical waves that disperse isotropically. 
Although there is a rather extensive literature devoted to the study of fully periodic elliptic operators, far less is known about partially periodic ones; we refer the reader to Section 9.2 of \cite{Kuchment} for some key references.
Most of these prior results appear to focus on the spectral properties of these operators, without a clear description of the associated dynamics beyond what is classically implied by the spectral type, as discussed below.
Since wave transport in periodic media is relatively well-understood, it is therefore natural to ask whether there is a parallel theory in waveguide settings.

In this paper, we are concerned with Schr\"{o}dinger operators on $\bbR^{n+m}$ with potentials that are compactly supported in the first $n$ coordinate directions, and periodic in the last $m$.
Such operators may produce \emph{surface states} \cite{black2022scattering,DaviesSimon}, which may be thought of as Bloch waves that are localized near the support of the potential. 
A famous result of Asch-Knauf \cite{AschKnauf} shows that solutions to the Schr\"{o}dinger equation on $\bbR^{d}$ with a fully periodic potential undergo ballistic transport.
One therefore expects that the surface states should evolve ballistically under the Schr\"{o}dinger evolution, but only along the periodic directions.
This phenomenon, referred to below as \emph{directional ballistic transport}, is the object of this paper.

The transport properties of fully periodic operators are a consequence of their Floquet theory.
Thus, to study partially periodic operators, we develop an analytic Floquet theory that parallels the classical fully periodic one.
In the latter setting, the Floquet transform conjugates a periodic elliptic operator on $\bbR^d$ to a direct integral of analytically varying operators, each acting on $L^2(\bbT^d)$, where $\bbT^d$ is the $d$-dimensional torus.
Because the torus is compact, $H(k)$ forms an analytic family of Fredholm operators, so one may use the theory of Fredholm determinants to parameterize their eigenvectors analytically \cite{wilcox1978theory}.
If $H$ is instead only partially periodic, there is still a (partial) Floquet transform, but the fiber operators $H(k)$ act on a cylinder, which is no longer compact, and therefore the operators $H(k)$ are not Fredholm.
The main challenge in the partially periodic setting is this lack of compactness, which makes the analytic variation of surface states more subtle than that of Bloch waves.
In particular, there may be eigenvalues embedded in the essential spectrum of each fiber operator, whose variation in $k$ is then a delicate matter.
We note that this is an issue particular to the open or ``leaky" waveguide setting in which there is no boundary condition in the transverse directions.
Many of the aforementioned prior works on waveguides pose the entire problem on a cylinder or a slab with boundary conditions, allowing one to bypass these considerations. 

Works in open waveguide settings, such as \cite{filonovKlopp,hoang2014absence,filonov}, have sought to rule out the possibility of eigenvalues or singular continuous spectrum under varying decay assumptions on the potential.
Though one can infer a weak dynamical corollary from such results via the RAGE theorem \cite{AmreinGeorgescu, Ruelle}, this is unsatisfactory because it cannot give any quantitative information about the propagation of surface states, let alone account for its directional character.
The strategy pursued in these references involves studying the variation of the resolvent in $k$ and appears to attain an analytic control of the surface states that is too weak to show ballistic transport; see Section~\ref{sec:prior} for a more thorough discussion.

Beyond that, with regard to dynamical results, we mention our previous work \cite{BDMY} in which directional ballistic transport results are proven on $\bbZ^2$. As we explain in Section~\ref{sec:prior}, the methods used there do not apply to the higher-dimensional continuum setting considered here.

To develop an analytic Floquet theory strong enough to prove ballistic transport, we instead appeal to a tool from the theory of elliptic boundary value problems: the Dirichlet-to-Neumann map.
This allows us to reformulate the eigenvalue problem for $H(k)$ as a Fredholm problem on a compact domain.
A similar reformulation of a related spectral problem on $\bbR^2$ (corresponding to $n=m=1$ below) appeared previously in work of Fliss \cite{fliss2013dirichlet}.
There, the focus was on establishing such an equivalence for eigenvalues outside of the essential spectrum of the fiber operators, with an eye toward numerical applications.  In contrast, our work considers all energies.

\subsection{Main results}
For $n\geq 2$, $m\geq 1$, label the first $n$ coordinates of $\bbR^{n+m}$ via $x$ and the last $m$ via $y$. Let $B_1$ be the unit ball in $\bbR^n$ and let $\e_1,\ldots\e_{n+m}$ be the standard basis vectors.

We consider a real-valued potential $V(x,y)\in L^\infty(\bbR^{n+m})$ satisfying:
\begin{align}
\label{eq:VasumpX}
    \supp V(x,y)\subset B_1\times \bbR^m,
\end{align}
and
\begin{align}
\label{eq:VasumpY}
        V((x,y)+\e_i)=V(x,y),\, \,\,\forall (x,y)\in\bbR^{n+m},\,n+1\leq i \leq m+n.
\end{align}

Denote by $H$ be the Schr\"{o}dinger operator on $L^2(\bbR^{n+m})$ associated to $V$, i.e.,
\begin{align*}
    &H=H_0+V,\\
    &H_0=-\Delta,
\end{align*}
which is self-adjoint with domain the Sobolev space $H^2(\bbR^{n+m})$.

The space of surface states is given in terms of the partial Floquet transform $U$, defined precisely in Section~\ref{sec:preliminaries}.
This transform conjugates $H$ to a direct integral over the Brillouin zone $\calB=[0,2\pi)^m$, so that
\begin{align*}
    UHU^*= (2\pi)^{-m}\int\limits_{\calB}^\oplus H(k)\:dk.
\end{align*}
Here, each $H(k)$ is a self-adjoint operator acting on $L^2(\bbR^n\times \bbT^m)$.
In terms of this decomposition, the space of surface states is given by,
\begin{align*}
    \calH_{\textrm{sur}}:=U^{*}\left((2\pi)^{-m}\int\limits_{\calB}^\oplus\calH_{\mathrm{pp}}(k)\:dk\right),
\end{align*}
where $\calH_{\textrm{pp}}(k)$ denotes the pure point subspace of $H(k)$.

For an operator $A$, we denote by $A_H(t)$ the Heisenberg-evolved operator $e^{itH}Ae^{-itH}$.
Let $X$ and $Y$ be the partial position operators
\begin{align*}
    X\psi(x,y)=(x_1,\ldots,x_n)\psi(x,y),\;Y\psi(x,y)=(y_1,\ldots,y_m)\psi(x,y),
\end{align*}
with the natural domains, and let $Q=(X,Y)$. We make the following definition:
\begin{definition}
    We say that $\psi\in D(Q)$ exhibits \emph{directional ballistic transport} under the evolution of $H$ if 
    \begin{align*}
        \lim_{t\rightarrow\infty}\frac{1}{t}X_H(t)\psi=0
    \end{align*}
    and
    \begin{align*}
        \lim_{t\rightarrow\infty}\frac{1}{t}Y_H(t)\psi
    \end{align*}
    exists and is non-zero.
\end{definition}

Our main theorem is then as follows:
\begin{theorem}
\label{thm:transport}
    Let $H$ be a Schr\"{o}dinger operator with potential $V$ satisfying \eqref{eq:VasumpX} and \eqref{eq:VasumpY}. Then every state $\psi \in D(Q)\cap H^2(\bbR^{n+m})\cap \calH_{\mathrm{sur}}\setminus\{0\}$ exhibits directional ballistic transport under the evolution of $H$.
\end{theorem}

From the aforementioned work of Asch-Knauf \cite{AschKnauf}, ballistic transport for periodic Schr\"{o}dinger operators is a consequence of the analytic variation of Bloch states with respect to the quasimomentum $k$.
Analogously, in our setting, Theorem~\ref{thm:transport} follows as a consequence of the analyticity result we prove for partially periodic operators. In fact,  proving the following theorem represents the main effort of this paper. 

\begin{theorem}
\label{thm:analyticity}
    Let $H$ be as in Theorem~\ref{thm:transport}. Then there are a countable collection of open sets $U_l\subset\calB $, functions $\lambda_{l,i}:U_l\rightarrow \bbR$, and finite rank projectors $k\in\calB\mapsto \pi_{l,i}(k):L^2(\bbR^n\times \bbT^m)\rightarrow L^2(\bbR^n\times \bbT^m)$ such that
    \begin{enumerate}
        \item\label{item:thm1} We may write
        \begin{align*}
            U H\vert_{\calH_{\mathrm{sur}}} U^*=\sum_{l,i} \int\limits_{U_l}^\oplus\lambda_{l,i}(k)\pi_{l,i}(k)\:dk.
        \end{align*}
        \item\label{item:thm2}  Each $\lambda_{l,i}$ and $\pi_{l,i}$ is analytic as a function of $k$ on the complement of a closed set of measure zero.
        \item\label{item:thm3} The set of $k$ for which $\nabla_k \lambda_{l,i}(k)=0$ is measure zero. 
    \end{enumerate}
\end{theorem}

Note that this Theorem may be vacuous if the surface subspace is empty.
It is easy, however, to construct examples of potentials that satisfy \eqref{eq:VasumpX} and \eqref{eq:VasumpY} that induce surface states, for instance, by choosing $V$ to be sufficiently negative.

Beyond transport, the analytic structure of fully periodic operators underpins their spectral and dynamical properties and accounts for much of the richness of their theory.
Thus, this theorem may prove useful in subsequent analyses of partially periodic operators.
For instance, Theorem~\ref{thm:analyticity} immediately recovers a special case of a result of Filonov and Klopp \cite{filonovKlopp} on the spectral type of $H$:
\begin{corollary}
\label{cor:ac}
    Let $H$ be as in Theorem~\ref{thm:transport}. Then its spectrum is purely absolutely continuous.
\end{corollary}

\subsection{Related work}\label{sec:prior}
Most prior work on partially periodic operators, at least in the ``open" setting, has sought to characterize their spectral type.
For $V(x,y)$ that is periodic in $y$, Filonov-Klopp \cite{filonovKlopp} showed that if $V$ decays superexponentially in $x$, then $H$ has absolutely continuous spectrum, a result they later extended to a similar class of Maxwell operators \cite{filonov2005absolutely}.
Under the weaker condition $|V(x,y)|<C|x|^{-1-\eps}$, Filonov \cite{filonov} showed the absence of eigenvalues for both Schr\"{o}dinger and Maxwell operators.
We also mention a work of Hoang-Radosz \cite{hoang2014absence}, in which it was shown that if a periodic potential on $\bbR^2$ is perturbed by a periodic line defect, then the resulting Maxwell or Schr\"{o}dinger operator has no eigenvalues.
This result is notable compared to the others because it accommodates a periodic background.
In all of these works, the basic idea is to study the resolvent of $H(k)$ as an analytic function of $k$ and the energy $E$ on an appropriate weighted Sobolev space.
One then shows that the poles of the resolvent are non-constant.
Because the poles could also correspond to resonances, such an argument does not seem sufficient to establish the analytic variation of eigenvalues.

With regard to the dynamics of surface states, at least in the continuum setting, we are aware only of the work of Davies-Simon \cite{DaviesSimon}, who proved that they obey a weak form of confinement. In particular, their Theorem~6.1 shows that for $n=1$, $\calH_{\textrm{sur}}$ may be equivalently characterized as those states that satisfy
\begin{align*}
    \lim_{a\rightarrow\infty}\sup_t\int_{|x|>a}|(e^{-itH}\psi)(x)|^2\:dx=0.
\end{align*}
In fact, it is known more generally from \cite{black2022scattering} that potentials that are merely concentrated in the $x$ directions without any periodicity assumption in $y$ generate a surface subspace that obeys a similar, though even weaker, confinement condition.
For ballistic transport of surface states, the only available result seems to be our previous work \cite{BDMY}.
There we showed that for an analogous class of Schr\"{o}dinger operators on $\bbZ^2$, the surface states exhibit directional ballistic transport.
Compared to the proofs of the above theorems, this is significantly easier to establish because the fiber operators act on a discrete cylinder with one extended direction, thus allowing the use of transfer matrices.

There is a rather large literature devoted to the study of quantum transport phenomena for other classes of Schr\"{o}dinger operators with continuous spectrum, which we make no attempt to comprehensively review, instead refering the reader to \cite{damanik2024ballistic} for an overview.
Here, we simply mention that proving transport for higher dimensional operators seems to be a hard problem.
In fact, even a perturbative result for short-range potentials on $\bbZ^d$ was only obtained quite recently \cite{damanik2025ballistic}.
Beyond this and the classic result of Asch-Knauf \cite{AschKnauf} on periodic operators (and its discrete analogue obtained by Fillman \cite{fillman2021ballistic}), the only other ballistic transport results in higher dimensions seem to be due to Karpeshina et al. \cite{Stolzetall, karpeshina2021ballistic}. 
In these works, ballistic lower bounds for the Abel mean of the position operator are obtained for certain quasi-periodic and limit-periodic operators on $\bbR^2$, and for generic quasi-periodic operators on $\bbR^d$, for $d\geq 2$, respectively.
From this perspective, our Theorem~\ref{thm:transport} is interesting because it establishes a strong form of ballistic transport for another class of higher-dimensional Schr\"{o}dinger operators. 

\subsection{Overview of the proof}
The Dirichlet-to-Neumann map (DtN), also referred to as the Poincar\'{e}-Steklov operator, is a classical object in the theory of boundary value problems; see, e.g., \cite{sylvester1990dirichlet} for a variety of applications to inverse problems.
Given an elliptic operator on a domain (or, more generally, manifold) with, say, smooth boundary and suitable boundary data, it solves the Dirichlet problem and then returns the trace of the solution's normal derivative back onto the boundary.
To study the eigenvalue problem for $H(k)$, we will divide the cylinder $\bbR^n \times \bbT^m$ into two domains: a compact interior and its complement. Thus, the operator is $H(k)-E$, for some $k\in \calB$ and $E\in\bbC$, and the manifold is one of these two domains.
For certain boundary data, we compute both the interior and exterior DtN and form the \emph{two-sided} DtN$, \Lambda(k,E)$, via their difference.
The basic observation behind the proof of Theorem~\ref{thm:analyticity} is that there is a correspondence between (boundary traces of) eigenfunctions of $H(k)$ of energy $E$ and the kernel of $\Lambda(k,E)$ - at least when the space of boundary data is suitably restricted.
The utility of this formulation is that, for the Laplacian on a bounded domain in $\bbR^d$, it is well-known that the corresponding DtN is elliptic \cite{mclean,taylor2010partial}.
Typically, this is proven via the method of \emph{layer potentials}, and by adapting it to our setting, we are able to show that this also holds for the two-sided DtN $\Lambda$, introduced above; see Section~\ref{sec:Fredholm}.
Thus, we are able to convert an eigenvalue problem for a non-elliptic operator into a Fredholm\footnote{Actually, because of the restriction on the data, $\Lambda$ is only left semi-Fredholm, but this weaker property suffices.} problem, albeit one that depends nonlinearly on the energy and quasimomentum.

Making this correspondence precise requires some effort because one must be careful to only select boundary data that yields $H^2$ solutions at infinity.
This is because the space of admissible data varies with the values of $E$ and $k$, for reasons that have to do with the spectral theory of the twisted Laplacian, $H_0(k)$, on the cylinder $\bbR^n\times \bbT^m$.
By decomposing in Fourier modes in the $y$ variables, one can think of this operator as a sum of Laplacians on $\bbR^n$, each shifted by a value depending on $k$ and the mode $j$; see Section~\ref{subsec:exterior}.
Given some energy $E$, one then expects to be able to solve the boundary value problem at infinity only on those modes for which $E$ is not in the spectrum of the shifted Laplacian.
The space of admissible data, denoted by $\calV(k,E)$ below, thus varies as $(k,E)\in \calB\times \bbR$ is moved across certain hypersurfaces.
On these hypersurfaces themselves, corresponding to the infimum of the spectrum on certain modes, the description of the relevant subspace is particularly involved, as we explain below.

These issues stem from the exterior problem, but we also briefly mention a difficulty related to the interior problem, namely the Dirichlet spectrum.
Indeed, the interior Dirichlet problem for each $k\in \calB$ is only guaranteed to be uniquely solvable for energies $E$ outside of the Dirichlet spectrum of $H(k)$ on $B_R\times \bbT^m$.
The appearance of this set is completely artificial in that it depends on the choice of cutoff radius $R$, which is used to define the compact interior.
Unsurprisingly, then, one can vary this parameter to ensure that any fixed $E$ does not lie in the spectrum of $H(k)$ for all $k$ in some open set.
This is carried out in Proposition~\ref{pr:RVariation}.

With the Fredholm property of the DtN in hand, one may obtain the analytic variation of surface states by using the theory of analytic varieties.
Arguments of this form originate in the work of  Wilcox \cite{wilcox1978theory} and have been extended in several directions by Kuchment and other authors; see \cite{kuchment1993floquet} and \cite{Kuchment}.
To complete the proof of  Theorem~\ref{thm:analyticity}, it remains to show that the functions $\lambda_{l,i}$ are non-constant almost everywhere.
The analogous statement in the fully periodic setting follows from an argument due to Thomas \cite{thomas}, which has served as the basis for most proofs of the continuity of the spectrum of periodic operators.
Typically, this involves obtaining lower bounds on the resolvent of $H_0(k)$ as the quasimomentum $k$ is complexified to have a large imaginary part.
Here too, we give a similar argument in Section~\ref{sec:thomas} by analyzing the free DtN (that is, when $V\equiv 0$). 
It is only in this section that we use the explicit form of this operator given in terms of Bessel functions.

Finally, given Theorem~\ref{thm:analyticity}, the directional ballistic transport result Theorem~\ref{thm:transport} follows from a similar implication in our previous work \cite{BDMY}.
This part of the proof transfers quite readily to the setting considered here, so we give only a sketch.

\subsection*{Outline of the paper} This paper is organized as follows: in Section~\ref{existenceSection}, we define the Floquet transform and introduce various spaces used in the analysis. In Section~\ref{Section:DtNConstruction}, we construct the DtN and establish the correspondence between its kernel and eigenfunctions. In Section~\ref{sec:Fredholm} we prove that the DtN is left semi-Fredholm. In Section~\ref{sec:theoremProofs} we obtain the proofs of the main theorems. Finally, in Appendix~\ref{sec:FreeComputation} we collect some results on boundary value problems on $\bbR^n$.
\section{Preliminaries}\label{existenceSection}
\subsection{Floquet Theory }
\label{sec:preliminaries}

Let $\bbT^m=\bbR^m/ \bbZ^m$ denote the torus and let $H^2(\bbR^n\times \bbT^m)$ be the usual Sobolev space.
Let $\calB:=[0,2\pi)^m$ be the Brioullin zone.
Define for $f\in L^2(\bbR^{n+m})$ and $k\in \calB$ the \emph{partial Floquet transform}
\begin{align*}
    (Uf)(x,y,k)=\sum\limits_{j\in \bbZ^m} e^{i\langle k, y+ j\rangle }f(x,y+j).
\end{align*}
 The standard properties of the Floquet transform, as recorded, for instance, in Section 4 of \cite{Kuchment}, extend readily to the partial Floquet transform:
\begin{proposition}\label{floquetPr}
    The partial Floquet transform $U$ has the following properties:
    \begin{enumerate}
        \item $U$ is a unitary map
        \begin{align*}
            U:L^2(\bbR^{n}\times \bbR^m)\rightarrow (2\pi)^{-m}\int\limits_{\calB}^\oplus L^2(\bbR^n\times \bbT^m)\:dk.
        \end{align*}
        \item We have the unitary equivalence
    	\begin{align*}
    		UHU^* =(2\pi)^{-m}\int\limits_{\calB}^\oplus H(k)\:dk,
    	\end{align*}
        where $H(k)$ is the self-adjoint operator
        \begin{align}\label{eq:H(k)Def}
        \begin{split}
            &H(k)=H_0(k)+V\\
            &H_0(k):=-\Delta-2ik\cdot \nabla_y+k^2
            \end{split}
        \end{align}
        with domain $H^2(\bbR^n\times \bbT)$.
    \end{enumerate}
\end{proposition}
Here we have used the shorthand that $k^2=\sum_{j=1}^mk_j^2$.

For later use, we will also need to consider $H_0(k)$ and $H(k)$ on the complexified Brillouin zone,
\begin{align*}
    \tilde{\calB}:=\{k\in \bbC^m\mid \Re k \in\calB\},
\end{align*}
where the real part is taken entrywise. The expressions in \eqref{eq:H(k)Def}
 then extend analytically to $\tilde{\calB}$. As in the fully periodic setting, $k_i\mapsto H_0(k)$ and $k_i\mapsto H(k)$ are analytic families for each index $i$ individually.

\subsection{Domain division and function spaces}
We divide the space $\bbR^n\times \bbT^m$ into two domains:
\begin{align*}
    \Omega^-_R=B_R\times \bbT^m,\quad
    \Omega^+_R=(\overline{B_R})^c\times \bbT^m,\quad
    \partial \Omega_R=\partial B_R\times \bbT^m,
\end{align*}
for $B_R$ the ball of radius $R$ in $\bbR^n$, $R>1$. 
We will often suppress the dependence of these sets on $R$ as we only need to vary it in Section~\ref{sec:analytic}. 

We let $H^2(\Omega^\pm)$ be the Sobolev spaces on $\Omega^\pm$ with norms
\begin{align*}
    &\|\varphi\|_{H^{2}(\Omega^\pm)}^2=\sum_{|\alpha|\leq 2}\|\partial^\alpha \varphi\|^2_{L^2(\Omega^\pm)}.
\end{align*}
We will also need the scale of Sobolev spaces $H^{s}(\partial\Omega)$, $s\geq 0$, whose norm is most conveniently expressed via the spectral decomposition of the Laplacian on $\partial\Omega$. 
To this end, first let
\begin{align*}
    \xi_j(y)=e^{i 2\pi  \braket{j,y}}
\end{align*}
for any $j\in \bbZ^m$ be the $j$-th Fourier mode in the $y$ directions. Now recall that the Laplacian on the unit sphere $\bbS^{n-1}$ has a complete set of orthonormal eigenfunctions (the spherical harmonics), usually denoted $Y_{\ell,m}(\omega)$, indexed by $\ell \in \bbN_0,m\in \calM_\ell$,\footnote{Note that we use the index $m$ here, since this is the common convention,  though it should not be confused with the dimension $m$.} where $\calM_\ell$ is some finite collection of indices with size depending on $n$. Each $Y_{\ell,m}$ corresponds to eigenvalue of $-\ell(\ell+n-2)$. It follows that (up to suitable normalization) $\{Y_{\ell,m}(x/R)\xi_j(y)\}_{\ell\geq 0,m\in \calM_\ell,j\in\bbZ^m}$ forms an orthonormal basis of $L^2(\partial\Omega)$ and for $f\in L^2(\partial\Omega)$, we let $f_{\ell,m,j}$ be its coefficients with respect to this basis. Then we define
\begin{align*}
    \|f\|_{H^s(\partial\Omega)}^2=\sum_{j\in\bbZ^m}\sum_{\ell\geq 0}\sum_{m\in\calM_\ell}(1+\ell^2+\|j\|^2)^s|f_{\ell,m,j}|^2.
\end{align*}

For $j\in\bbZ^m$, we also introduce the closed subspace $H_j^s(\partial \Omega)\subset H^s(\partial \Omega)$ that consists of functions of the form $\psi(x)\xi_j(y)$ and define $H_j^2(\Omega^\pm)$ analogously. We also let
\begin{align*}
    \calL_{\ell,j}=\vspan_{m\in\calM_\ell}\{f(r)Y_{\ell,m}(\omega)\xi_j(y)\},
\end{align*}
which lies in $H^s(\partial\Omega)$ for any $s\geq 0$.
Below, we will need the traces from the interior and exterior domains as well as the corresponding extension operators. The traces $T^\pm$ act on smooth functions via restriction, and the extension operators $\Ex^\pm$ are their right inverses. The following is standard, see, e.g., \cite[Ch. 4]{taylor1996partial}.
\begin{proposition}
    The interior and exterior trace operators $T^-$ and $T^+$ extend as bounded operators
    \begin{align*}
        &T^\pm:H^2(\Omega^\pm)\to H^{3/2}(\partial\Omega).
    \end{align*}
    Moreover, they each have a right inverse $\Ex^\pm$, which is a bounded operator
\begin{align*}
        \Ex^\pm:H^{3/2}(\partial\Omega)\to H^{2}(\Omega^\pm).
    \end{align*}
\end{proposition}

Below, we will also use $T$ to denote the trace from $\bbR^n\times \bbT^m$ to $\partial\Omega$.
\section{Construction of the Dirichlet-to-Neumann map}\label{Section:DtNConstruction}
In this section, we associate to most pairs $(k, E)\in \tilde{\calB}\times \bbR$ a Dirichlet-to-Neumann map (DtN), denoted $\Lambda(k, E)$, and show that its kernel consists of boundary traces of eigenfunctions of $H(k)$ with energy $E$. 
To do this, we first examine the solvability of the interior and exterior Dirichlet problems for the operator $H(k)-E$.

\subsection{The interior boundary value problem}
To begin, we let $H^-_{0,D}(k,R)$ be the operator $H_0(k)$ acting on $L^2(\Omega^-_R)$ (the interior set), with domain
\begin{align*}
    \calD(H_{0,D}^-)=\{u\in H^2(\Omega^-)\mid T^-u=0\},
\end{align*}
i.e., with Dirichlet boundary condition. Similarly, we define $H_{D}^-(k)$, to be the operator $H(k)$ acting on the same domain. We suppress the variable $R$, and record the following elementary proposition. 
\begin{proposition}\label{pr:DirichletSA}
    For each $k\in \bbR^m$, the operators $H^-_{0,D}(k)$ and $H^-_D(k)$ are self-adjoint and are closed for any $k\in\bbC^m$. Moreover, they have discrete spectra that accumulate only at infinity.
\end{proposition}
\begin{proof}
    For $k\in \bbC^m$, $H_{0,D}^-(k)$  is a relatively bounded perturbation of $H_{0,D}^-(0)$. It follows that $H_{0,D}^-(k)$ is closed for all $k\in \bbC^m$, and self-adjoint for $k\in \bbR^m$ by the Kato-Rellich theorem. Since $H_{0,D}^-(k)$ has a compact resolvent by Rellich-Kondrachov, so does $H_{D}^-(k)$, guaranteeing the spectra are discrete, with accumulation only at infinity.
\end{proof}
This allows us to define the resolvent
\begin{align*}
    R_D^-(k,E)=(H^-_D(k)-E)^{-1},
\end{align*}
as a bounded operator from $L^2(\Omega^-)\rightarrow H^2(\Omega^-)$, whenever $E\not\in \sigma(H_D^-(k))$. We use this to give the following result on the interior Dirichlet problem.
\begin{proposition}\label{pr:innerBVP}
    Let $E\in\bbC\setminus \sigma(H_D^-(k))$ for some $k\in \tilde{\calB}$. Then for all $f\in H^{3/2}(\partial \Omega)$ there exists a unique $u\in H^2(\Omega^-)$ so that    
    \begin{align}\label{eq:innerBVP}
    \begin{cases}
        &H(k)u=Eu\,\text{in }\Omega^-\\
        &T^-u=f.
        \end{cases}
    \end{align}
    The solution map $f\mapsto u$ denoted $D^-(k,E)$ is bounded from $H^{3/2}(\partial\Omega)\to H^2(\Omega^-)$ and may be given by
    \begin{align}\label{eq:D-form}
        D^-(k,E)=(\Id-R_D^-(k,E)(H^-(k)-E))\Ex^-.
    \end{align}
\end{proposition}
\begin{proof}
    That 
    \begin{align*}
        u=(\Id-R_D^-(k,E)(H(k)-E))\Ex^-f
    \end{align*}
    solves \eqref{eq:innerBVP} is readily verified. Uniqueness follows from the fact that the solution to
    \begin{align*}
        \begin{cases}
                &H^-(k)u=Eu\\
                &T^-u=0.
        \end{cases}
    \end{align*}
    is solved uniquely by $R_{0,D}^-(k,E)0=0$, for $E\not \in \sigma(H_D^-(k))$.
    The boundedness of 
    \begin{align*}
        D^-:H^{3/2}(\partial \Omega)\to H^2(\Omega^-)
    \end{align*}
    follows immediately from the expression \eqref{eq:D-form}.
\end{proof}
\subsection{The exterior boundary value problem }
\label{subsec:exterior}
Next, we turn to the exterior problem, that is, the problem on $\Omega^+$. 
The solvability of the exterior boundary value problem requires that the boundary data lie in a subspace of $H^{3/2}(\partial\Omega)$ that depends on $k$ and $E$. To understand this, note that $H_0(k)$ acts on each subspace $H^2_j(\Omega^{+})$ as
\begin{align}
\label{eq:jAction}
    -\Delta_x+(k+2\pi j)^2.
\end{align}
In fact, if we let $H_{0,D}^+(k)$ be $H_0(k)$ on $L^2(\Omega^+)$ with domain,
\begin{align*}
    \calD(H_{0,D}^+(k))=\{u\in H^2(\Omega^+)\mid T^+u=0\}
\end{align*}
then 
\begin{align*}
    \sigma(H_{0,D}^+(k)\vert_{ H_j(\Omega^+)})=[0,\infty)+(k+2\pi j)^2.
\end{align*}
This suggests that, for any $k$ and $E$, the exterior Dirichlet problem for $H_0(k)$ at energy $E$ should only be solvable if the data has no components in those $H_j^{3/2}(\partial \Omega)$ for which 
\begin{align*}
     E_j(k):=E-(k+2\pi j)^2
\end{align*}
is in $[0,\infty)$.
Accordingly, for any $(k,E)\in \calB\times \bbR$ we define the subspace of $H^{3/2}(\partial\Omega)$
\begin{align*}
    \calV_-(k,E)=\bigoplus_{\{j\in \bbZ^{m} \mid E_j(k)<0\} }H^{3/2}_j(\partial \Omega),
\end{align*}
and we let $\tilde{R}^+_D(k,E)$ be the resolvent of $H_{0,D}^+$ restricted to this subspace at energy $E$. As for the interior Dirichlet resolvent, $\tilde{R}^+_D(k,E)$ is bounded from $L^2(\Omega^+)\rightarrow H^2(\Omega^+)$. Note that as a function of $k$ and $E$, the space $\calV_-(k,E)$ is constant on the complement of the hypersurfaces $E_j(k)=0$.

In Section~\ref{sec:thomas}, we will need to complexify $k$, so for $k\in \tilde{\calB}$, we let
\begin{align*}
    \calV_-(k,E)=\bigoplus_{\{j\in \bbZ^{m} \mid E_j(\Re k)<0\} }H^{3/2}_j(\partial \Omega),
\end{align*}
where the real part is taken entrywise. We note that with this definition $\calV(k,E)$ does not change as $k$ acquires an imaginary part. In this case,
\begin{align*}
    E_j(k)=E_j(\Re k) +(\Im k)^2 -2 i \Im k\cdot (\Re k+2\pi j),
\end{align*}
so $\tilde{R}^+_D(k,E)$ is still well-defined if the final term is non-zero.

It turns out that if $E_j(k)=0$, the solvability of the exterior boundary value problem in $H^{2}_j(\Omega^+)$ is more subtle and depends on which angular momentum sectors the data lies in, at least in low dimensions. Thus, we define
\begin{align*}
    &\calV_0(k,E):=\bigoplus_{\{j\in\bbZ^{m} \mid E_j(k)=0 \} }\bigoplus_{\ell\geq \ell_0(m)}\calL_{\ell,j},\\
    &\ell_0(n):=\begin{cases}
        2&n=2\\
        1&n=3,4\\
        0&n\geq 5.
    \end{cases}
\end{align*}
This subspace is constant along the hypersurfaces $E_j(k)=0$, at least where no two hypersurfaces intersect.

For a given $(k,E)$, the subspace
\begin{align}
\label{eq:calVDef}
    \calV(k,E):=\calV_-(k,E)\oplus \calV_0(k,E)
\end{align}
is exactly the space of data for which the exterior boundary value problem is solvable:
\begin{proposition}\label{pr:outerBVP}
    Fix $E\in\bbR$ and $k\in \tilde{\calB}$. Then for $f\in H^{3 /2}(\partial \Omega) $ there exists $u\in H^2(\Omega^+)$ so that    
    \begin{align}
    \label{eq:outerBVP}
        \begin{cases}
            &H(k)u=Eu\,\text{in }\Omega^+\\
    	&T^+u=f
            \end{cases}
    \end{align}
    if and only if $f \in \calV(k,E)$. Moreover, the solution, when it exists, is unique.
    The solution map $f\mapsto u$ denoted $D^+(k,E)$ is bounded from $\calV(k,E)\to H^2(\Omega^+)$ and is given explicitly by
    \begin{align}
    \label{eq:outerBVPSol}
    \begin{split}
        u(x,y)&=(\Id - \tilde{R}_D^+(k,E)(H(k)-E))\Ex^+ \tilde{f}\\
    &\quad+  \sum_{\{j\in \bbZ^m\mid E_j(k)=0\}}\sum_{\ell\geq \ell_0}\sum_{m\in\calM_l}f_{\ell,m,j}\left(\frac{r}{R}\right)^{2-n-\ell}Y_{\ell,m}(\omega)\xi_j(y),
    \end{split}
    \end{align}
    for $\tilde{f}$ the component of $f$ in $\calV_-$ and $r=|x|$.
\end{proposition}
For the proof, we need the following Lemma on the exterior Dirichlet problem in $\bbR^n$. Its proof is given in Appendix~\ref{sec:FreeComputation}.
\begin{restatable*}{lemma}{exteriorZeroRn}
\label{lem:exteriorZeroRn}
    Let $f\in H^{3/2}(\partial B_R)$ be given by
    \begin{align*}
        f=\sum_{\ell\geq 0}\sum_{m\in\calM_\ell}f_{\ell,m}Y_{\ell,m}. 
    \end{align*}
    Then there exists $u\in H^2(B_R^c)$ solving
    \begin{align}
    \label{eq:outerDeltaBVP}
    \begin{cases}
        -\Delta u=0\\
            T^+ u=f
    \end{cases}
    \end{align}
    if and only if $f_{\ell,m}=0$ for all $\ell<\ell_0(n)$.
    Moreover, the solution is unique, is given by
    \begin{align*}
        u(r\omega)=\sum_{\ell\geq \ell_0}\sum_{m\in \calM_\ell}f_{\ell,m}\left( \frac{r}{R} \right)^{2-n-\ell}Y_{\ell,m}(\omega),
    \end{align*}
    for $r\in \bbR^+,\omega \in \bbS^{n-1} $, 
    and satisfies
    \begin{align}
    \label{eq:zeroSolMappingProperty}
        \|u\|_{H^2(B_R^c)}\leq C \|f\|_{H^{3/2}(\partial B_R)}
    \end{align}
    for some dimensional constant $C>0$.
\end{restatable*}

\begin{proof}[Proof of Proposition~\ref{pr:outerBVP}]
   It is easy to check that if $f\in \calV(k,E)$, then $u$ as defined in \eqref{eq:outerBVPSol} is a solution to \eqref{eq:outerBVP} and by Lemma~\ref{lem:exteriorZeroRn} and the argument in Proposition~\ref{pr:innerBVP} it must be unique.
   The fact that
   \begin{align*}
       \|u\|_{H^2(\Omega^+)}\leq C\|f\|_{H^{3/2}(\partial\Omega)}
   \end{align*}
   is a consequence of \eqref{eq:zeroSolMappingProperty} and the mapping properties of $\tilde{R}_D^+(k,E)$.

   To see that there can only be a solution if $f\in\calV$, observe that
    the boundary value problem \eqref{eq:outerBVP} is equivalent to
    \begin{align}
    \label{eq:outerBVPj}
    \begin{cases}
        &(-\Delta_x - E_j(k))u_j=0\,\text{in }B_R^c\subset \bbR^n\\
        &T^+u_j=f_j
    \end{cases}
    \end{align}
    for all $j\in \bbZ^m$.
    It follows from a theorem of Rellich that if $E_j(k)\in(0,\infty)$, then the only $L^2$ solution to \eqref{eq:outerBVPj} is $0$, see \cite[Lem. 9.8]{mclean} for the relevant version. On the other hand, if $E_j(k)=0$, then $f_{\ell,j,m}=0$ for all $\ell<\ell_0(n)$ by Lemma~\ref{lem:exteriorZeroRn}. Since these are exactly the conditions that $f\in \calV$, we are done.
\end{proof}

\subsection{The two-sided Dirichlet-to-Neumann map}
We may now define the two-sided Dirichlet-to-Neumann map.
\begin{definition}
    Let $\calS=\{(k,E)\in \tilde{\calB}\times \bbR\mid E\in \sigma(H_{0,D}^-(k))\}$.
    For any $(k,E)\in (\tilde{B}\times\bbR)\setminus \calS$, we define the \emph{two-sided Dirichlet-to-Neumann} map $\Lambda:\calV(k,E)\rightarrow H^{1/2}(\partial\Omega)$
    \begin{align*}
        &\Lambda(k,E):=\Lambda^-(k,E)-\Lambda^+(k,E)\\
        &\Lambda^\pm(k,E):=T^\pm\partial_rD^\pm(k,E).
    \end{align*}
    Here, $\partial_r$ denotes the radial derivative with respect to the $x$ variables.
\end{definition}
To show that this map indeed encodes the existence of eigenfunctions, we recall the following version of Green's identity, see for instance \cite[Lem. 4.4]{mclean}:
\begin{lemma}\label{GreensId}
    Suppose that $u,v\in H^1(\Omega^\pm)$. If $\Delta u$ and $\Delta v$ are in $L^2(\Omega^\pm)$, then
    \begin{align*}
        \left<u,\Delta v\right>_{\Omega^\pm}=\left<\Delta u,v \right>_{\Omega^{\pm}}\pm \left<T^{\pm}\partial_r u,T^\pm v \right>_{\partial\Omega}\mp \left<T^\pm u,T^\pm\partial_rv \right>_{\partial\Omega}.
    \end{align*}
	Similarly, if $(H(k)-E)u$ and $(H(\overline{k})-E)v$ are in $L^2(\Omega^\pm)$, then 
    \begin{align}\label{GreenIdEq}
        \left<u,(H(k)-E)v\right>_{\Omega^\pm}=\left<(H(k)-E)u,v \right>_{\Omega^{\pm}}\pm \left<T^{\pm}\partial_r u,T^\pm v \right>_{\partial\Omega}\mp \left<T^\pm u,T^\pm\partial_rv \right>_{\partial\Omega}.
    \end{align}
    In particular, if $\varphi\in C^\infty_c(\bbR^n\times \bbT^m)$ is a test function, then
    \begin{align}\label{eq:Greens}
    	\left<D^\pm f,(H^{\pm}(k)-E)\varphi \right>_{\Omega^{\pm}}=\pm \left<T^{\pm}\partial_rD^\pm f,T\varphi \right>_{\partial\Omega}\mp\left<f,T\partial_r\varphi\right>_{\partial \Omega}.
    \end{align}
\end{lemma}
With this, we may establish the key property of $\Lambda$, which is that its kernel corresponds to eigenvectors of $H(k)$ with eigenvalue $E$.

\begin{proposition}\label{pr:IsoProp}
    Fix $(k,E)\in(\tilde{\calB}\times \bbR)\setminus \calS$,
    the map
    \begin{align*}
	    \calV(k,E)\ni f\mapsto (D^-(k,E)+D^+(k,E))f\in H^2(\bbR^{n}\times \bbT^{m})
    \end{align*}
    is an isomorphism between $\ker\Lambda(k,E)$ and the eigenspace of $H(k)$ corresponding to $E$.
\end{proposition}
\begin{proof}
	First, we verify that if $\Lambda(k,E)f=0$, then $u=(D^-+D^+)f$ is an eigenvector of $H(k)$ with eigenvalue $E$. If $\psi\in C_0^\infty(\bbR^{n}\times \bbT^{m})$, then by \eqref{eq:Greens}
    \begin{align}
    	\begin{split}\label{eq:pmGreens}
    	\left<u,(H(k)-E)\psi \right>&=\left<D^+f,(H(k)-E)\psi \right>_{\Omega^+}+\left<D^-f,(H(k)-E)\psi \right>_{\Omega^-}\\
    	&=\left<T^+\partial_rD^\pm f,T\psi \right>_{\partial\Omega}-\left<T^-\partial_rD^\pm f,T\psi \right>_{\partial\Omega},
    	\end{split}
    \end{align}
    which vanishes since $f\in \ker\Lambda$. Thus, $(H(k)-E)u$ is zero as a distribution, so we need only check that it is in $H^2(\bbR^{n}\times \bbT^m)$. This is clear from the definition of $D^\pm$ on $(\partial\Omega)^c$. In a neighborhood of $\partial\Omega$,  $u$ solves 
    \begin{align*}
    	(H_0(k)-E)u=0	
    \end{align*}
    in the sense of distributions so it is in fact smooth there by local elliptic regularity \cite[Prop. 3.9.1]{taylor1996partial}.\par
    This confirms that $u$ is an eigenfunction. To see every eigenfunction arises this way, suppose that $(H(k)-E)u=0$. Then because $u\in H^2(\bbR^{n}\times \bbT^m)$, we have that $T^+u=T^-u$, which we denote $f\in H^{3/2}(\partial \Omega) $. By the uniqueness statements of Propositions \ref{pr:innerBVP} and \ref{pr:outerBVP}, we see that
    \begin{align*}
    	u=(D^++D^-)f
    \end{align*}
    and furthermore from \eqref{eq:pmGreens} it follows that $\Lambda f=0$. Finally, we note that because $D^+T^+f$ is $L^2$ at infinity, $f$ must lie in $\calV(k,E)$ by Proposition~\ref{pr:outerBVP}. Therefore, $\varphi\mapsto f$ is an inverse, so we are done.
\end{proof}

\section{Fredholm property of the Dirichlet-to-Neumann map}
\label{sec:Fredholm}
Having established the connection between $\Lambda$ and eigenfunctions of $H(k)$, we now turn to its Fredholm properties.
Because $\Lambda$ is defined only on the subspace $\calV$, it will be only left semi-Fredholm, i.e., it has a finite-dimensional kernel and closed range, but its cokernel may be infinite dimensional.
To see this, we use the method of layer potentials.
Typically, on a finite volume domain, the single-layer potential is defined via the Green's function of $\Delta$ on the full space, but we prefer to instead use the resolvent kernel at energy $-1$ to avoid some technicalities at infinity.
To this end, let $G=(-\Delta_{\bbR^{n}\times \bbT^{m}} +1)^{-1}$ and let $G(x,y)$ be its kernel.
Recall from \cite[(7.53)]{teschl} that for $z \not\in [0,\infty)$, the kernel of $(-\Delta_{\bbR^{d}}-z)^{-1}$ for an arbitrary $d\geq 2$ is given by
\begin{align}
\label{eq:RdGreensZ}
    G_{\bbR^{d}}(w;z)=\frac{1}{2\pi}\left( \frac{\sqrt{-z}  }{2\pi |w|}\right) ^{\frac{d}{2}-1}K_{\frac{d}{2}-1}(\sqrt{-z} |w|),
\end{align}
where $K_{\nu}$ is the modified Bessel function of the second kind, and the root is always the one with a positive real part.
We will also let $G_{\bbR^d}(w;0)$ be the standard Green's function of $-\Delta$, that is,
\begin{align}
\label{eq:RdGreens0}
    G_{\bbR^d}(w;0)=\frac{\Gamma(d/2)}{2\pi^{d/2}(d-2)}|w|^{2-d}.
\end{align}
By decomposing in Fourier modes in $y$, we have that 
\begin{align*}
    G(x,y)=\sum_{j\in\bbZ^m}G_{\bbR^n}(x;-1-\|2\pi j\|^2)\xi_j(y).
\end{align*}
Now let $T^*:H^{-1/2}(\partial\Omega)\rightarrow H^{-1}(\bbR^n\times \bbT^m)$ be the adjoint (with respect to the $L^2$ inner product) of the trace operator acting from $H^{1}(\bbR^n\times \bbT^m)\rightarrow H^{1/2}(\partial\Omega)$.
For $f\in H^{-1/2}(\partial\Omega)$, we define the \emph{single-layer operator} $\SL$ via
\begin{align*}
    \SL f:=GT^{*} f.
\end{align*}

It is clear from the definition that $(-\Delta+1)\SL f$ is  $0$ (in the sense of distributions) on $\Omega^{\pm}$, as $f$ is supported on $\partial \Omega$.
Since $G$ is smooth for $(x,y)\neq 0$, when $f\in L^1(\partial \Omega)$ we have the integral representation
\begin{align}
\label{eq:SLRep}
    (\SL f)(x,y)=\int_{\partial\Omega}G(x-x',y-y')f(x',y')\:dx'dy',
\end{align}
for $(x,y)\in \bbR^n\times \bbT^m\setminus \partial\Omega$.

The following behavior of the resolvent kernel is well-known and is essentially a consequence of the fact that  $G$ is a pseudo-differential operator of order  $-2$.
\begin{proposition}
    In a neighborhood of $(0,0)$, the kernel $G(x,y)$ obeys 
    \begin{align}
        \label{eq:closeToR^nm}
        G(x,y)=G_{\bbR^{n+m}}(x,y;0)+\calR(x,y),
    \end{align}
    where $\calR$ satisfies
    \begin{align}
    \label{eq:kernelBoundR}
        |\partial_x^\alpha\partial_y^\beta \calR(x,y)|\leq C_{\alpha,\beta}d((x,y),0)^{3-n-m-|\alpha|-|\beta|},
    \end{align}\
    for all multi-indices $\alpha,\beta$.
    In particular, we have that
    \begin{align}
        \label{eq:kernelBound}
        |\partial_x^\alpha\partial_y^\beta G(x,y)|\leq C_{\alpha,\beta}d((x,y),0)^{2-n-m-|\alpha|-|\beta|},
    \end{align}
    near $0$.
\end{proposition}
\begin{proof}
    Let $\chi(x,y)$	be a smooth function that is identically $1$ on a neighborhood of $(0,0)\in  \bbR^{n}\times \bbT^{m}$ and compactly supported inside a chart of this manifold. 
    We may then regard $G\cdot \chi$ as a distribution on $\bbR^{n+m}$. 
    Because the cylinder $\bbR^{n}\times \bbT^{m}$ is flat, we have that
    \begin{align*}
        (-\Delta_{\bbR^{n+m}} +1)(G\cdot\chi)=\delta_0+\nabla G\cdot \nabla \chi+G\cdot \Delta\chi.
    \end{align*}
    The latter two terms are smooth because $G$ is smooth away from $(0,0)$ and $\chi$ is $1$ in a neighborhood of $(0,0)$.
    Therefore, 
    \begin{align*}
            (-\Delta_{\bbR^{n+m}} +1)(G(x,y)\cdot\chi(x,y) -G_{\bbR^{n+m}}(x,y;1))
    \end{align*}
    is smooth so that $G(x,y)\cdot\chi(x,y) -G_{\bbR^{n+m}}(x,y;1)$ is smooth by elliptic regularity. It thus suffices to show that $G_{\bbR^{n+m}}(x,y;1)-G_{\bbR^{n+m}}(x,y;0)$ is singular of order $3-n-m$.
    
    To prove this, we need the standard Bessel function asymptotics as $w\to 0$ from, e.g., \cite[\S 10.31]{AS}.
    They are
    \begin{align}
        \label{eq:BesselAsym1}
        K_{\nu}(w)=\frac{1}{2}(\frac{1}{2}w)^{-\nu}(\Gamma(\nu)+O(w^{2}))
        +(-1)^{\nu+1}(\frac{1}{2}w)^{\nu}\log(\frac{1}{2}w)(1+O(w^{2})),
    \end{align}
    when $\nu$ is a positive integer and 
    \begin{align}
        \label{eq:BesselAsym2}
        K_{\nu/2}(w)=\frac{1}{2}\Gamma(\nu /2)(\frac{1}{2}w)^{-\nu /2}(1+O(w^2)),
    \end{align}
    when $\nu$ is an odd positive integer, with all $O(w^{2})$ terms  analytic.
    Inserting these asymptotics into \eqref{eq:RdGreensZ}, we find that if $n+m$ is even, then
    \begin{align*}
        G_{\bbR^{n+m}}(x,y;1)&=\frac{\Gamma(\frac{n+m}{2}-1)}{4\pi^{(n+m)/2}}|(x,y)|^{2-n-m}(1+O(|(x,y)|^2)\\
        &\quad+(-1)^{(n+m)/2-1)}\frac{1}{2\pi}\log(\frac{1}{2}|(x,y)|)(1+O(|x,y|^2)),
    \end{align*}
    whereas if $n+m$ is odd, we may omit the second term.
    Comparing to \eqref{eq:RdGreens0}, we obtain \eqref{eq:closeToR^nm}, and \eqref{eq:kernelBound}  follows immediately by differentiating.
\end{proof}

The bound \eqref{eq:kernelBound} implies that $G(x-x',y-y')$ is uniformly integrable on $\partial \Omega$ as ${(x,y)\to (x',y')}$ so that $(T^+\SL f)(x,y)$ and $(T^-\SL f)(x,y)$  agree and are given by the expression \eqref{eq:SLRep} for $(x,y)\in \partial\Omega$.\par
To proceed, we will need the following mapping properties.
\begin{proposition}
\label{pr:SLProperties}
    The operator $\SL$ is bounded from $H^{-1/ 2}(\partial\Omega)\to H^{1}(\bbR^{n}\times \bbT^{m})$.
    Moreover, $T^{\pm}\SL$ is bounded from $H^{1 /2}(\partial \Omega)\to H^{2}(\Omega^{\pm})$.
\end{proposition}
\begin{proof}
    We will treat $\Omega^{+}$, with $\Omega^{-}$ being similar.
    By duality, $T^{*}$ is bounded $H^{-1 /2}(\partial\Omega)\to H^{-1}(\bbR^{n}\times \bbT^{m})$. Since $G$ is bounded $H^{-1}(\bbR^{n}\times \bbT^{m})\rightarrow H^{1}(\bbR^n\times \bbT^m)$, we obtain the mapping property for $\SL$.
    
    From \eqref{eq:kernelBound}, we know that on the $n+m-1$ dimensional manifold $\partial\Omega$, $G(x-x',y-y')$ defines the Schwartz kernel of a pseudo-differential operator of order  $-1$ because the singularity is of order $1-(n+m-1)$, see \cite[Prop. 2.4]{taylor2010partial}. Therefore, $T^{\pm}\SL$ maps $H^{1 /2}(\partial\Omega)\to H^{3 /2}(\partial\Omega)$, as claimed.
\end{proof}
Unlike $\SL f$, which is continuous, the normal derivative of $\SL f$ exhibits a discontinuity across $\partial\Omega.$ The analogous jump relation for the Green's function of the Laplacian is well-known.
\begin{proposition}[Jump relation]
\label{pr:Jump}
    We have the identity
    \begin{align*}
    	T^{\pm}\partial_r \SL=\mp \frac{1}{2}\Id + N,
    \end{align*}
    where $N$ is the integral operator
    \begin{align*}
    Nf(x,y)=\int_{\partial \Omega}\partial_rG(x-x',y-y')f(x',y')\:dx'dy'
    \end{align*}
    for $(x,y)\in \partial\Omega$.
\end{proposition}
\begin{proof}
    From the integral representation of $\SL$, we have that 
    \begin{align*}
    	\partial_r \SL f&=\int_{\partial\Omega} \partial_rG(x-x',y-y')f(x',y')\:dx'dy'\\
    			&=\int_{\partial\Omega}\partial_rG_{\bbR^{n+m}}(x-x',y-y')f(x',y')\:dx'dy'+\int_{\partial\Omega}\partial_r\calR(x-x',y-y')f(x',y')\:dx'dy',
    \end{align*}
    for $(x,y)\not\in\partial\Omega$.
    Due to \eqref{eq:kernelBoundR}, $\partial_r \calR$ is uniformly integrable as  $(x',y')\to (x,y)\in\partial\Omega$.
    Thus, the function defined a priori for $(x,y)\notin \partial\Omega$ by the second integral extends to $\partial\Omega$, and therefore has a two-sided trace given by the same expression for $(x,y)\in \partial \Omega$.
    On the other hand, we have the classic jump relation for the single-layer potential \cite[Prop. 7.11.3]{taylor2010partial}
    \begin{align*}
    &T^{\pm}\int_{\partial U}\partial_r G_{\bbR^n\times \bbR^m}(x-x',y-y')f(x',y')\:dx'dy'=\\
    &\quad\mp \frac{1}{2}f+\int_{U}\partial_r G_{\bbR^n\times \bbR^m}(x-x',y-y')f(x',x')\:dx'dy',
    \end{align*}
    valid for $U$ any domain in $\bbR^{n+m}$ with smooth boundary.
    By passing to a partition of unity subordinate to charts on $\partial\Omega$ and applying this identity locally, the proof is now complete.
\end{proof}

\begin{proposition}\label{pr:SLLambda}
    We have the following identities:
    \begin{align}
        &T^{+}\SL \Lambda^{+}=  -\frac{1}{2}\Id_{\calV} +N+\calK_1\label{eq:regularizer1}\\
        &T^{-}\SL \Lambda^{-}= \frac{1}{2}\Id_{\calV} +N+\calK_2\label{eq:regularizer2}
    \end{align}
    where
    \begin{align*}
        &\calK_1=T^+G (2ik\cdot \nabla-k^2+E+1)D^+\\
        &\calK_2=T^-G (2ik\cdot \nabla-k^2-V+E+1)D^-.
    \end{align*}
\end{proposition}
\begin{proof}
    Let $f\in \calV(k,E)$ and $g\in L^2(\partial \Omega)$.  
    For \eqref{eq:regularizer1}, we use that $\SL=G T^*$ and take adjoints to write
    \begin{align*}
    	\left<T^{+}\SL \Lambda^{+}f,g \right>_{\partial\Omega}&=\left<\Lambda^{+}f,T\SL g\right>_{\partial \Omega}\\
        &=\left<T^{+}\partial_ru^{+},T\SL g \right>_{\partial\Omega},
    \end{align*}
    where we have set $u^{+}=D^{+}f$.
    Since $u^{+}$ and $\SL g$ are both in $H^1(\Omega^{+})$ with $\Delta u^{+}$ and $\Delta \SL g$ in $L^2(\Omega^+)$, we are free to apply Green's identity (Lemma~\ref{GreensId}) to find that
    \begin{align*}
        \left<T^{+}\partial_ru^{+},T^{+}\SL g \right>_{\partial\Omega}=
        \left< f,T^{+}\partial_r\SL g\right>_{\partial \Omega}+\left<u^{+},\Delta\SL g \right>_{\Omega^+}-\left<\Delta u^{+},\SL g \right>_{\Omega^+}.
    \end{align*}
    The first term may be rewritten via Proposition~\ref{pr:Jump} as
    \begin{align*}
        \left< f,T^{+}\partial_r\SL g\right>_{\partial \Omega}=\left<(-\frac{1}{2}+N)f, g\right>_{\partial\Omega}.
    \end{align*}
    Using that $\Delta \SL g=\SL g$ inside $\Omega^+$ we may also write
    \begin{align*}
        \left<u^+,\Delta \SL g \right>_{\Omega^{+}}=\left<u^+,\SL g \right>_{\Omega^{+}}
        =\left<T^{+}G u^+,g \right>_{\partial\Omega}.
    \end{align*}
    Similarly, because $(H_0(k)-E)u^+=0$, we have that
    \begin{align*}
    	\left<\Delta u^{+},\SL g \right>_{\Omega^{+}}
    =\left<T^{+}G (-2ik\cdot \nabla+(k^2-E))  u^+,g \right>_{\partial\Omega}.
    \end{align*}
    Combining these three identities and applying duality yields \eqref{eq:regularizer1}.
    
    The derivation of \eqref{eq:regularizer2} is the same except that we obtain an extra term in $\calK_2$ from instead using
    \begin{align*}
    	(H_0(k)+V-E)D^{-}f=0.
    \end{align*}
\end{proof}
\begin{corollary}
\label{cor:Fredholm}
For any $(k,E)\in(\tilde{\calB}\times \bbR)\setminus \calS$, the full DtN $\Lambda(k,E)$ is left semi-Fredholm. 
\end{corollary}
\begin{proof}
    By Proposition~\ref{pr:SLProperties}, we have that the operators $T^{\pm}\SL$ are bounded from $H^{1 /2}(\partial\Omega)\to \calV(k,E)$. 
    Moreover, by Proposition~\ref{pr:SLLambda}
    we have that
    \begin{align*}
    T\SL \Lambda = \Id_{\calV}-\calK_1+\calK_2.
    \end{align*}
    The operator $\Lambda$ will therefore have $T\SL$ as a left regularizer,  if we can show that $\calK_1$ and $\calK_2$ are compact.
    
    We have the following mapping properties for the composition defining $\calK_2$ :
    \begin{align*}
    H^{3/2}(\partial\Omega)\xrightarrow{D^{-}} H^{1}(\Omega^{-})\xrightarrow{-2ik\cdot\nabla -k^2-V+E+1} L^{2}(\Omega^{-})\xrightarrow{T^{-}G}H^{3 /2}(\partial\Omega).
    \end{align*}
    As we have seen, $D^{-}$ is bounded from $H^{3 /2}(\partial\Omega)\to H^{2}(\Omega^{-})$ and is therefore compact as an operator $H^{3 /2}(\partial\Omega)\to H^{1}(\Omega^{-})$. Thus, the first arrow is compact so that $\calK_2$ is as well. The same argument shows that $\calK_1$ is compact, so we are done.
\end{proof}

\section{Proofs of the main theorems}\label{sec:analytic}
\label{sec:theoremProofs}
In order to prove Theorem~\ref{thm:analyticity}, we first require some results that allow us to parameterize the kernel of analytic families of left semi-Fredholm operators.
Related results for families of Fredholm operators may be found in \cite{zaidenberg1975banach}.

\subsection{Analytical variation of the kernel of left semi-Fredholm operators}

\begin{proposition}\label{pr:simpleFredholmKernelAnalyticity}
Let $B_{\eps}\subset \bbT^{d}$ be an open convex subset and let
\begin{align*}
	k\in B_{\eps}\mapsto F(k):X\to Y
\end{align*}
be an analytic family of left semi-Fredholm operators between Hilbert spaces $X$ and $Y$.
Then there exists an open set of full measure $\mathcal{W}\subset B_\eps$, whose complement is analytic, and an analytic, with respect to $k\in \mathcal{W}$, operator-function projecting onto $\ker F(k)$.
\end{proposition}
\begin{proof}
    Consider a complex $\delta$-width (with small $\delta$) tubular neighborhood $U$ of $B_\epsilon$. It is geometrically, and thus also a holomorphically convex domain in $\bbC^n$. If the width $\delta$  is small, then the operator-function is still analytic and, due to the openness of the set of semi-Fredholm operators (see, e.g., \cite[Section VI.7]{Dunford} and \cite[Proposition XI.2.4]{Conway}), also semi-Fredholm.
    
    Let now $d_0\geq 0$ be the smallest dimension of $\ker F(k)$ in $U$, attained at a point $k_0$. Let us pick a complement $E_0$ to $\ker F(k_0)$. Consider the inverse operator $G$ to $F(k_0): E_0 \to \im F(k_0)$ and a bounded projector $P$ onto the range $\in F(k_0)$. Also, introduce the analytic operator-function
    $$B(k):=GPF(k).$$
    If the neighborhood $U$ is small enough, this is clearly an analytic Fredholm operator-function such that $\ker B(k)=\ker F(k)$ and thus $\dim\ker B(k)\geq d_0$. 
    
    Now one can refer to known results about Fredholm analytic operator-functions. Namely,
    let $A_0\subset U$ be the set where $\dim\,\ker F(k)\geq d_0+1$. According to Theorem 4.11 in \cite{zaidenberg1975banach}, the subset is analytic. Consider any non-zero analytic function $\phi$ from the ideal corresponding to $A_0$ and define the set $A\subset U$ as the set of all zeros of $\phi$. Clearly, $A_0\subset A$. The advantage of using $A$ is that it has co-dimension $1$ (unlike $A_0$, which could have larger co-dimension). Then $\mathcal{W}:=U\setminus A$ is holomorphically convex, as the complement of a principal analytic set in a domain of holomorphy. Moreover, for $k\in \mathcal{W}$, the kernel of $F(k)$ has the same dimension $d_0$. Thus, these kernels form an analytic $d_0$-dimensional sub-bundle in the ambient space, and thus an analytic in $\mathcal{W}$ projector exists locally (for instance, one can use the Riesz projector). Due to the holomorphic convexity of $\mathcal{W}$, such an analytic projector exists also globally (M. Shubin's theorem, see \cite{shubin1979holomorphic} or \cite[Theorem 3.11]{zaidenberg1975banach}, \cite[Theorem 1.6.13]{kuchment1993floquet}). 
    
    The intersection of the ``bad'' set $A$ with the real subspace has zero measure; see \cite[Lemma 5.22]{Kuchment} (not an original source) and references therein, in particular the nice derivation in \cite{mityagin2020zero}.
    
    This proves the theorem: indeed, one has existence of an analytic projector over a set $\mathcal{W}\bigcap\bbR^n$ of full measure.
\end{proof}

\begin{proposition}\label{pr:complicatedFredholmKernelAnalyticity}
For some $\eps>0$ and real numbers $a<b$, let 
\begin{align*}
	(k,E)\in B_{\eps}\times (a,b)\mapsto F(k,E): X\to Y
\end{align*}
be an analytic family of left semi-Fredholm operators between Hilbert spaces $X$ and $Y$.
Suppose that for each $k\in B_{\eps}$, the set
\begin{align*}
	\mathcal{D}_k:=\{E\in (a,b) \mid \dim \ker F(k,E)>0  \} 
\end{align*}
is discrete.
Then there exist countably many functions $\lambda_n(k):B_{\eps}\to (a,b)$ and families of projection operators 
\begin{align*}
k\in B_\eps\mapsto \pi_n(k):X\to X
\end{align*}
such that
\begin{enumerate}
	\item For each $n$, both $\lambda_n(k)$ and  $\pi_n(k)$ are defined and analytic almost everywhere in $B_{\eps}$.
	\item For all $k\in B_{\eps}$ 
		\begin{align*}
			\ker F(k,\lambda_n(k))=\Ran \pi_n(k).
		\end{align*}
	\item $\dim\ker F(k,E)>0$ if and only if $E=\lambda_n(k)$ for some $n$.
\end{enumerate}
\end{proposition}
\begin{proof}
	 \footnote{See also Theorem 5.5.1 in the forthcoming book \cite{Kuchment2025Periodic}.} 
Let us define the analytic variety
\begin{align*}
	\mathcal{D}:=\{k\in \mathbb{T}^m,E\in (a,b) \mid \dim \ker F(k,E)>0  \}.
\end{align*}
It is sufficient to consider one of its irreducible components $C$, which is automatically Stein. According to 
\cite[Theorem 1.6.14]{kuchment1993floquet}, the sheaf of the germs of analytic sections of $\ker\,(H(k)-\lambda I)$ over $C$ is coherent. Let the minimal dimension of this kernel over $C$ be equal to $d$. Then, due to the irreducibility of $C$, the dimension is the same outside of an analytic subset $S\subset C$. In a neighborhood of any point $(k_0,\lambda_0)\in C\setminus S$, there is a local holomorphic basis of sections $\phi_{J}(\lambda, k), j=1,\dots k$. Then, according to a classical approximation theorem (see, e.g. \cite[Theorem 7.2.7]{Hormander}), these functions can be in a (maybe smaller) neighborhood of $(k_0,\lambda_0)$ uniformly approximated by global sections $\psi_j$. Then, aside from a proper analytic subset $A$ of $B$, they stay the desired basis sections.
\end{proof}

Now recall that the subspace $\calV(k,E)$ defined in \eqref{eq:calVDef} is constant on any subset of $\calB\times \bbR$ where none of the $E_j(k)$ change sign. On such a neighborhood, $\Lambda$ is an analytic function of $(k,E)$ away from  $\calS$, the set of $(k,E)$ corresponding to the Dirichlet spectrum.
\begin{proposition}
    Let $U\subset \calB\times \bbR$ be a neighborhood on which $\calV(k,E)$ is constant. Then
    \begin{align*}
        (k,E)\in (U\times \bbR)\setminus \calS \mapsto \Lambda(k,E): \calV\rightarrow H^{1/2}(\partial\Omega)
    \end{align*}
    is an analytic family of left semi-Fredholm operators.
\end{proposition}
\begin{proof}
    The operators $H_{0,D}^\pm(k)$ are analytic families on $\calB$. Therefore, by \cite[Thm. XII.7]{RSVol4} their resolvents are analytic functions on $(U\times \bbR)\setminus \calS$. The claim then follows from the representations \eqref{eq:D-form} and \eqref{eq:outerBVPSol}.
\end{proof}

We now show that by varying $R$, we may avoid the set of $(k,E)$ for which $\Lambda(k,E)$ is not defined.
Let $\lambda_1(k,R)\leq \lambda_2(k,R)\leq\ldots$ be the eigenvalues of $H_{D}^-(k,R)$ with multiplicity.
\begin{proposition}
    \label{pr:RVariation}
    For any $(k_0,E_0)\in \calB\times \bbR$, there exists $R>1$ and a neighborhood $U\subset \bbT\times \bbR$ containing $(k_0,E_0)$ such that $E\not \in  \sigma(H_D(k,R))$ for all $(k,E)\in U$.
\end{proposition}
This will come as an easy consequence of the following lemma:
\begin{lemma}
	Fix $k\in \bbT^{*}$. For every $n\in \bbN$, the function $R\mapsto \lambda_n(k,R)$ is continuous and strictly monotone decreasing.
\end{lemma}
\begin{proof}
    We first establish the monotonicity. 
    For this, recall the min-max expression for eigenvalues with respect to form domains:
    \begin{align}
    \label{eq:min-max}
    	\lambda_n(R)&=\min_{\substack{\{\psi_1,\ldots,\psi_n\} \subset H^{1}_{0}(\partial\Omega_R)\\\|\psi_j\|=1}}\max_{\substack{\psi\in\vspan \{\psi_1,\ldots,\psi_j\} \\\|\psi\|=1}}\left<H_{D,R}(k)\psi,\psi \right>,
    \end{align}
    with the max achieved only by an eigenfunction.
    For $R'>R$, we have the natural inclusion  $H_0^1(\partial\Omega_R)\hookrightarrow H_0^1(\partial \Omega_{R'}) $, due to the Dirichlet boundary condition, via extension by zero, so we at least have that $\lambda_n(R)\geq \lambda_n(R')$.
    
    Now suppose that $\lambda_n(R)=\lambda_n(R')$ for some $R'>R$ and let $\psi_1,\ldots,\psi_n$ be $L^2$-normalized eigenfunctions for $\lambda_1(R),\ldots,\lambda_n(R)$. 
    Extending them by $0$, they achieve the min-max in \eqref{eq:min-max} for $H_{D,R'}(k)$ and therefore $\psi_n$ solves  $H(k)\psi_n=\lambda_n(R)\psi_n$ on $\Omega_{R'}^-$.
    In particular, it is zero on $\Omega_{R'}^-\setminus\Omega_R^-$ so that it must in fact be identically zero by unique continuation for Schr\"{o}dinger eigenfunctions, see, e.g., \cite[Thm. XIII.57]{RSVol4}. 
    This is a contradiction, as $\psi_n$ is a non trivial eigenfunction, so we must have that $\lambda_n(R')>\lambda_n(R)$, as claimed.
    
    To see the continuity, we use that the operator
    \begin{align*}
    \tilde{H}_{D}(k,R):=-R^{-2}\Delta_x -\Delta_y - 2i k \cdot \nabla_y+k^2+V(Rx,y),
    \end{align*}
    with domain $H_0^1(\Omega_1)\cap H^{2}(\Omega_1)$ has the same spectrum as $H_{D}(k,R)$.
    By Sobolev embedding, let $p>2$ be an exponent so that $H^{1}_0(\partial \Omega_1)\hookrightarrow L^{p}(\partial \Omega_1)$ and let $q$ solve $\frac{1}{q}+\frac{2}{p}=1$.
    Suppressing $k$, we observe that for any $\psi \in H^{1}_0(\partial \Omega_1)$
    \begin{align*}
        &\left|\left<\tilde{H}_{D}(R)\psi,\psi \right>-\left<\tilde{H}_{D}(R')\psi,\psi \right> \right|\\
        &\quad \leq \int_{ \Omega_{1}}(R^{-2}-(R^{'})^{-2})|\nabla_x \psi|^{2}+|V(Rx,y)-V(R'x,y)||\psi|^{2}(x,y)\:dxdy\\
        &\quad \leq (R^{-2}-(R^{'})^{-2})\|\psi\|_{H^{1}_0( \Omega_1)}^{2}
        +\|V(Rx,y)-V(R'x,y)\|_{L^{q}( \Omega_1)}\|\psi\|^{2}_{H^{1}_0( \Omega_1)},
    \end{align*}
    where in the last line we have used H\"{o}lder's inequality and Sobolev embedding.
    Since
    \begin{align*}
    	\|V(Rx,y)-V(R'x,y)\|_{L^{q}(\partial \Omega_1)}\to 0,\,R'\to R,
    \end{align*}
    for instance by taking a smooth approximant of $V$ in $L^{q}$, we have shown that
    \begin{align*}
    	\sup_{\psi \in H_0^{1}(\Omega_1),\|\psi\|_{H^1_0( \Omega_1)}=1}|\left< (\tilde{H}_{D}(R)-\tilde{H}_{D}(R'))\psi,\psi\right>|\to 0,\,R'\to R.
    \end{align*}
    By \cite[Thm. VIII.25]{RSVol1}, it follows that any sequence $\tilde{H}_{D}(R_j)$ with $R_j\to R$ is norm resolvent convergent. Therefore, by \cite[Thm. VIII.23]{RSVol1} the spectral projections to any interval $P_{(a,b)}(\tilde{H}_{D}(R_n))$ converge in norm to $P_{(a,b)}(\tilde{H}_{D}(R))$, from which the continuity of $\lambda_n(R)$ is a simple consequence. With this, the proof is complete.
\end{proof}
Proposition~\ref{pr:RVariation} now follows quickly:
\begin{proof}[Proof of Proposition~\ref{pr:RVariation}]
	For some $R>0$, let $\lambda_N(k_0,R)$ be the smallest eigenvalue of $H_{D}(k_0,R)$ that is greater than $E_0$. Because $\lambda_N(k_0,R)$ is continuous, if we increase $R$ by a sufficiently small amount,  it will still not be equal to $E_0$ and by monotonicity, the same will be true of all $\lambda_j(k_0,R)$ for $j<N$. Therefore, there is some $R'$ so that $E_0 \not\in \sigma(H_{D}(k_0,R))$. We may now find a neighborhood of $k_0$ on which $E_0\not\in \sigma(H_{D,R}(k))$ simply by standard analytic perturbation theory in $k$. Choosing a suitably small neighborhood of $E_0$, the proof is now complete.
\end{proof}
Now we apply these results to the DtN, yielding items \eqref{item:thm1} and \eqref{item:thm2} of Theorem~\ref{thm:analyticity}.
\begin{proposition}
\label{pr:analyticFamily}
	There exists a countable collection of open sets $U_l\subset \calB$ such that the following holds: for each $U_l$ there are finitely many families of finite rank projections
    \begin{align*}
        k\in U_l\mapsto \pi_{l,j}(k):H^2(\bbR^n\times \bbT^m)\rightarrow H^2(\bbR^n\times \bbT^m)
    \end{align*}
    and functions $\lambda_{l,j}(k):U_l\to \bbR$ such that
\begin{enumerate}
    \item The function $\lambda_{l,j}(k)$ and the family of operators $\pi_{l,j}(k)$ are analytic on the complement of a closed set of zero measure.
    \item For every $k\in U_l$ and every $j$, $\ker [H(k)-\lambda_{l,j}(k)]=\Ran \pi_{l,j}(k)$.
    \item If for some $(k_0,E_0)\in \calB \times \bbR$, $\ker [H(k_0)-E_0)]\neq 0$, then $E_0=\lambda_{l,j}(k_0)$ for some indices $l$ and $j$.
\end{enumerate}
\end{proposition}
\begin{proof}
We will show that for every $(k_0,E_0)\in \calB\times \bbR$ such that $\ker [H(k_0)-E_0]\neq 0$, there is some $U\subset \calB$ containing $k_0$ and $\lambda_j(k)$ and $\pi_j(k)$ satisfying the conclusions of the Proposition. By Proposition~\eqref{pr:RVariation}, we may vary $R$ so that $(k_0,E_0)\subset U\times (a,b)$ that is disjoint from $\calS$.

There are two cases. The first is when $(E_0)_j(k_0)\neq 0$ for all $j\in \bbZ^m$. In this case, we may restrict $U\times(a,b)$ so that $\calV(k,E)$ is constant on it. Then by Corollary~\ref{cor:Fredholm} and Proposition~\ref{pr:analyticFamily}, $\Lambda(k,E)$ is an analytic family of left semi-Fredholm operators on $U\times (a,b)$ so we may apply Proposition~\ref{pr:complicatedFredholmKernelAnalyticity} to produce functions $\lambda_j(k)$ and projectors $\pi_j(k)$ that characterize its kernel on $U\times (a,b)$ and are analytic off of a set of measure zero. The range of
\begin{align*}
    (D^{-}(k,\lambda_j(k))+D^{+}(k,\lambda_j(k)))\circ \pi_j(k)	
\end{align*}
for each index $j$ then varies analytically in $k$ almost everywhere. 
By Proposition~\ref{pr:IsoProp} the projector onto their range forms a complete set of eigenprojectors for $H(k)$ with eigenvalues in $E$ when $(k,E)\subset U\times(a,b)$. This completes the proof in this case.

For the second case, we suppose that $(E_0)_j(k_0)=0$ for some $j\in\bbZ^m$.
We may assume that $(E_0)_{j'}(k_0)\neq 0$ for any other $j'$ because the set of $k$ which this is not true is of measure $0$. Thus, we may restrict $U$ so that $\calV(k,E_0)=\calV(k_0,E_0)$ for all $k\in U$. 
We may now conclude by applying
Proposition~\ref{pr:simpleFredholmKernelAnalyticity} to the analytic family $k\in U\mapsto \Lambda(k,(E_0)_j(k))$ and arguing as in the first case.
\end{proof}


\subsection{Non-constancy of the energies}
\label{sec:thomas}
 Now we prove that none of the $\lambda_n(k)$ are locally constant. This shows that the spectrum of $H$ is purely absolutely continuous, thus proving Corollary~\ref{cor:ac}.

For this, we will analyze the free DtN, $\Lambda_0$, corresponding to $V\equiv 0$.
More precisely, we set
\begin{align*}
    &\Lambda_0(k,E):=\Lambda_0^-(k,E)-\Lambda^+(k,E)\\
    &\Lambda_0^-(k,E):= T^- \partial_r D_0^-(k,E),
\end{align*}
where $D_0^-(k,E)$ is the Dirichlet solution map for $H_0(k)-E$ on $\Omega^-$ (recall that $\Lambda^+$ has no potential by definition). 
Clearly $\Lambda_0^-(k,E)$  is well-defined for $E\not\in \sigma(H_{0,D}^-(k))$.

Based on an ODE computation in spherical harmonics given in Appendix~\ref{sec:FreeComputation}, we have the following representation for $\Lambda_0$ in terms of $P_\nu(w):=I_\nu(w)K_\nu(w)$, with $I_\nu$ the modified Bessel function of the first kind.

\begin{proposition}
    \label{pr:explicitLambda}
    Let $(k,E)\in \tilde{\calB}\times \bbR$. For any $j\in \bbZ^m$ such that $E_j(\Re k)<0$, we have that
    \begin{align*}
        (\Lambda_0(k,E) f)_{j}=-\sum_{\ell\geq 0}\sum_{m\in\calM_\ell}\frac{f_{\ell,m,j}}{RP_{n/2+\ell-1}(z_jR)}Y_{\ell,m}(\omega),
    \end{align*}
    where $z_j^2=-E_j(k)$, $\Re z>0$.
\end{proposition}
\begin{proof}
    By Lemma~\ref{lem:RnHelmholtz},
    \begin{align*}
        &(D_0^-f)_j=\sum_{\ell\geq 0}\sum_{m\in \calM_\ell}f_{\ell,m,j}\left(\frac{r}{R}\right)^{1-n/2}\frac{I_{n/2+\ell-1}(z_jr)}{I_{n/2+\ell-1}(z_jR)}\\
        &(D^+f)_j=\sum_{\ell\geq 0}\sum_{m\in \calM_\ell}f_{\ell,m,j}\left(\frac{r}{R}\right)^{1-n/2}\frac{K_{n/2+\ell-1}(z_jr)}{K_{n/2+\ell-1}(z_jR)}.
    \end{align*}
    Thus
    \begin{align*}
                &
                (\Lambda_0^- f)_j=\sum_{\ell\geq 0}\sum_{m\in\calM_\ell}\frac{f_{j,\ell,m}}{R}[(1-\frac{n}{2})+z_jR\frac{I_{n/2+\ell-1 }'(z_jR)}{I_{n/2 + \ell -1 }(z_jR)}]Y_{\ell,m}(\omega)\\
                &(\Lambda^+ f)_j =\sum_{\ell\geq 0}\sum_{m\in\calM_\ell}\frac{f_{j,\ell,m}}{R}[(1-\frac{n}{2})+z_jR\frac{K_{n/2 + \ell -1 }'(zR)}{K_{n/2 + \ell -1 }(z_jR)}]Y_{\ell,m}(\omega)
    \end{align*}
    so that
    \begin{align*}
        &(\Lambda_0 f)_j=\\
        &\quad-\sum_{\ell\geq 0}\sum_{m\in\calM_\ell}f_{\ell,m,j}z_j\frac{K_{n/2+\ell-1}(z_jR)I_{n/2+\ell-1}'(z_jR)- K_{n/2+\ell-1}'(z_jR)I_{n/2+\ell-1}(z_jR)}{P_{n/2+\ell-1}(z_j
        R)}Y_{\ell,m}(\omega),
    \end{align*}
    By \cite[(10.28.2)]{AS}, the Wronskian in the numerator evaluates to $1/(z_jR)$, from which the conclusion follows. 
\end{proof}

Next, we prove the following technical estimate:
\begin{proposition}\label{pr:PnuEst}
\label{pr:PNuEst}
    We have the following bound, valid for any $\nu \geq 0$ and $z\in\bbC$ with $\Re z>0$:
    \begin{align*}
        |P_\nu(z)|\leq \frac{C}{(|\Im z^2|)^{\frac{1}{3}}},
    \end{align*}
    for some universal constant $C>0$. 
\end{proposition}
\begin{proof}
    We start with the following identity, that holds for $\Re{z}>0,\Re(\nu) >-1 $, see \cite[(6.535)]{Zwillinger2007Table} :
    \begin{align*}
        P_\nu(z)=\int\limits_0^\infty \frac{x}{x^2+z^2}[J_\nu(x)]^2\, dx.
    \end{align*}
    Next we use the following bound from \cite{landau2000bessel}, that holds for $x\in \bbR,\nu\geq 0 $
    \begin{align*}
        |J_\nu(x)|<c|x|^{-\frac{1}{3}},
    \end{align*}
    and a uniform constant $c>0$.
    Combining these and setting $a=\Im z^2$, we get that
    \begin{align*}
        |P_\nu(z)|&\leq c\int\limits_0^\infty \frac{x^{\frac{1}{3}}}{|x^2+z^2|}\: dx\\
        &\leq c\int\limits_0^\infty \frac{x^{\frac{1}{3}}}{\sqrt{x^4+a^2}}\:dx\\
        &=c|a|^{-1}\int\limits_0^\infty \frac{x^{\frac{1}{3}}}{\sqrt{(x^2/|a|)^2+1}}\:dx\\
        &= c|a|^{-1/3}\int\limits_0^\infty \frac{u^{\frac{1}{3}}}{\sqrt{u^4+1}}\:du,
    \end{align*}
    from which the conclusion follows immediately.
\end{proof}

To show that there are no flat bands in the spectral variety of $H$, we give an argument inspired by the classical method of Thomas \cite{thomas} that relies on complexifying the quasimomentum.
\begin{lemma}\label{lem:nonconstancy}
    For any $E_0\in\bbR$, the set of $k\in\calB$ for which $E_0$ is an eigenvalue of $H(k)$ is a set of measure zero. 
\end{lemma}
\begin{proof}
    Let $E_0 \in \bbR$, we will show that for any $k_0\in [0,2\pi)^m$ with $|(k_0)_1|>0$, there is some interval $0\in I$ so that
    \begin{align}\label{SetOfEval}
         \{k\in k_0+I\mathbf{e_1}\mid E_0 \text{ is an eigenvalue of }H(k)\}
    \end{align}
    is of Lebesgue measure zero.
    The conclusion then follows from Fubini's theorem and a covering argument.\par
    Note that each set $\{(k,E)\in \calB\times \bbR\mid E_j(k)=0\}$ is a paraboloid so the measure of 
    \begin{align*}
        \{k_0\in \calB\mid (E_0)_j(k_0)=0\,\text{for some }j\in\bbZ^m\}
    \end{align*}
    is zero.
    Therefore, we may assume that $(E_0)_j(k_0)$ is non-zero for each $j$.
    Moreover, by Proposition~\ref{pr:RVariation} we may choose $R$ such that $E_0$ is not in $\sigma(H_{0,D}^-(k))$ or $\sigma(H_D^-(k))$ for any $k\in k_0+I\e_1$, where $I$ is some sufficiently small closed interval.
    This guarantees that $\calV(k,E_0)=\calV(k_0,E_0)$ for all $k\in k_0+I\e_1$.
    We further assume that $I$ is small enough so that all $k\in k_0 +I\e_1$ satisfy $|k_1|>\delta$ and $|k_1-2\pi|>\delta$ for some $\delta>0$. 

   Now let $\bbV$ be the vertical strip in the complex plane
    \begin{align*}
        \bbV=\{z\in \bbC\mid \Re z\in I\}.
    \end{align*}
    Note that by the analytic Fredholm theorem \cite[\S 8.1]{yafaev1992mathematical},
     the set of $z\in \bbV$ for which 
    \begin{align*}
        E_0\in \sigma(H_{0,D}^-(k_0+z))\cup \sigma(H_D^-(k_0+z))
    \end{align*}
     is discrete. Furthermore, the conditions on $k_1$ guarantee that $\tilde{R}_D^+(k,E_0)$ is well-defined for all $k\in k_0+I\e_1$.
     
     Thus, the operators given by
    \begin{align*}
        &\tilde{\Lambda}_0(z):=\Lambda_0(k_0+z,E_0)\\
        &\tilde{\Lambda}(z):=\Lambda(k_0+z,E_0)
    \end{align*}
    are an analytic family of Fredholm operators on $\bbV$ up to a discrete set.

    From Proposition~\ref{pr:explicitLambda}, we have that
    \begin{align*}
        (\tilde{\Lambda}_0(it)f)_j=-\sum_{\ell\geq 0}\sum_{m\in\calM_\ell}\frac{f_{\ell,m,j}}{RP_{n/2+\ell-1}(z_j(t)R)}Y_{\ell,m}(\omega),
    \end{align*}
    where we used
    \begin{align*}
        (z_j^2(t))&=-(E_0)_j(k_0+it\e_1)=-E_0+(k_0+2\pi j+it\e_1)^2\\
        &=-(E_0)_j(k_0)-t^2 + 2it((k_0)_1+2\pi j_1)
    \end{align*}
    and $\Re z_j(t)>0$. Since $|(k_0)_1|>\delta$ and $|(k_0)_1-2\pi |>\delta$, it follows that $|\Im z_j^2(t)|>C|t|$ uniformly in $j$, so that by Proposition~\ref{pr:PNuEst}, we have that
    \begin{align*}
    \frac{1}{|P_{n/2+\ell-1}(z_j(t)R)|}\geq C|t|^{1/3},
    \end{align*}
    uniformly in $\ell$ and $j$.
    Thus, 
\begin{align*}
    \|\tilde{\Lambda}_0(it)f\|_{H^{1/2}(\partial \Omega)}\geq C|t|^{1/3}\|f\|_{H^{3/2}(\partial \Omega)}
\end{align*}
    for any $f\in \calV$.  

    Referring to \eqref{eq:D-form}, we write
    \begin{align*}
        D^-(k,E)-D_0^-(k,E)&=[-R_D\cdot(H(k)-E)+R^-_{0,D}\cdot(H_0(k)-E)]\Ex^-\\
        &=\left[(R^-_{0,D}-R^-_D)((H(k)-E))-R^-_{0,D}V\right]\Ex^-\\
        &=R^-_{0,D}V(R^-_D\cdot(H(k)-E)-\Id)\Ex^-,
    \end{align*}
    where we have used the resolvent identity in the last line.
    Now we use that $R_{0,D}(k_0+it,E_0)$ is bounded from $L^2(\Omega^-)\rightarrow H^2(\partial\Omega^-)$ independently of $t$
    to see that
    \begin{align*}
        \|\tilde{\Lambda}(it)-\tilde{\Lambda}_0(it)\|_{\calV\rightarrow H^{1/2}(\partial \Omega)}\leq C
    \end{align*}
    for all $t$ sufficiently large and therefore $\tilde{\Lambda}(it)$ is invertible for all $t$ sufficiently large.

    To conclude, we note that by the analytic Fredholm theorem, $\tilde{\Lambda}(z)$ is invertible for all but a discrete set of points in $\bbV$. Therefore, the set of $k\in k_0+I\e_1$ for which $\Lambda(k)$ has a non-trivial kernel is of measure zero, which is what we wanted to show.
\end{proof}

The proof of Theorem~\ref{thm:analyticity} is now immediate:
    \begin{proof}[Proof of Theorem~\ref{thm:analyticity}]
    Parts \eqref{item:thm1} and \eqref{item:thm2} follow from Proposition~\ref{pr:analyticFamily} whereas \eqref{item:thm3} follows from Lemma~\ref{lem:nonconstancy}.    
\end{proof}

To obtain Theorem~\ref{thm:transport}, we appeal to some arguments given in \cite{BDMY}. 
Though the setting considered there is different, the proofs are similar enough that we merely highlight the differences.
\begin{proof}[Proof of Theorem~\ref{thm:transport}] 
    Fix $\psi \in D(Q)\cap H^2(\bbR^{n+m})\cap \calH_{\mathrm{sur}}\setminus\{0\}$.
    First, must show that $\lim_{t\rightarrow\infty}\frac{1}{t}Y_H(t)\psi$ exists and is non-zero.
    From a standard limiting argument (see, e.g., the proof of \cite[Thm. 2.3]{AschKnauf}), it suffices to establish the theorem for $\psi$ of the form
        \begin{align}
            \calU\psi=\left(\sum_{l=1}^N\int\limits_{U_l}^\oplus\sum_{n=1}^M\pi_{l,n}(k)\:\frac{dk}{\abs{\calB}}\right)\calU\psi
        \end{align}
        for some $N,M>0$. 
    By \cite{RadinSimon}, we have the identity
        \begin{align*}
            Y_H(T)\psi =Y_H(0)\psi+2\int\limits_{0}^T P^y_H(t)\psi \:dt,
        \end{align*}
    where $P^y$ is the momentum operator $-i(\partial_{y_1},\ldots,\partial_{y_m})$, for $\psi$ in this domain.
    It therefore suffices to show that 
    \begin{align*}
        \lim_{T\rightarrow\infty}\frac{1}{T}\int\limits_{0}^T P^y_H(t)\psi \:dt
    \end{align*}
    exists and is non-zero. The proof of Proposition 4.1 of \cite{BDMY} now applies verbatim to our setting, so that we may conclude, as the input from  Theorem~\ref{thm:analyticity} is identical in either setting.
    
    Now, we must show that $\lim_{t\rightarrow\infty}\frac{1}{t}X(t)\psi=0$.
    This is simply Proposition 4.2 of \cite{BDMY}. Note that though the proof is written for a Schr\"{o}dinger operator on $\bbZ^2$, it applies without alteration to the present setting. 
\end{proof}

\appendix

\section{Boundary value problems on $\bbR^n$}
\label{sec:FreeComputation}

In this Appendix, we collect some simple facts about the Dirichlet problem for the Laplace and Helmholtz equations on the ball $B_R^-$ and its exterior $B_R^+=\overline{B_R^-}^c$ in $\bbR^n$. Although these results are likely standard, we were unable to locate a suitable reference; therefore, we include proofs for the reader's convenience.
\exteriorZeroRn

\begin{proof}
We decompose
\begin{align*}
    u(r\omega)=\sum_{\ell\geq 0}\sum_{m\in \calM_\ell}u_{\ell,m}(r)Y_{\ell,m}(\omega).
\end{align*}
By writing $\Delta$ in spherical coordinates, we find that each $u_{\ell,m}$ satisfies the ODE 
\begin{align*}
r^2u''_{\ell,m}(r)+(n-1)ru'(r)-\ell(\ell+n-2)u(r)=0,
\end{align*}
which has a basis of solutions given by $\{r^{\ell},r^{2-n-\ell}\} $ (unless $n=2$ and $\ell=0$ in which case a basis of solutions is given by $\{1,\log r\}$). Recall that
\begin{align*}
\ell_0(n):=\begin{cases}
    2&n=2\\
    1&n=3,4\\
    0&n\geq 5,
\end{cases}
\end{align*} 
and observe that $r^{2-n-\ell}\in L^2([R,\infty),r^{n-1}\:dr)$ if and only if $\ell\geq \ell_{0}$, which forces $f_{\ell,m}=0$ for $\ell<\ell_0$.
We therefore have that
\begin{align*}
    u(r\omega)=\sum_{\ell\geq \ell_0}\sum_{m\in \calM_\ell}f_{\ell,m}\left( \frac{r}{R} \right)^{2-n-\ell}Y_{\ell,m}(\omega),
\end{align*}
solves~\eqref{eq:outerDeltaBVP}, so it remains to estimate $\|u\|_{H^{2}(B_R^{+})}$.
By scaling, it suffices to take $R=1$.

Using the orthogonality of the $Y_{\ell,m}$, we integrate to find that
\begin{align*}
    \|u\|_{L^2(B_R^c)}^2&\leq C\sum_{\ell\geq\ell_0}\sum_{m\in\calM_\ell}^{} \ell^{-1}\|f_{\ell,m}\|^2\\
    &\leq C\|f\|_{L^2(\partial B_R)}^{2}.
\end{align*}
Writing the gradient in spherical coordinates, we obtain
\begin{align*}
    \nabla u&=\frac{\partial u}{\partial r}\hat{r}+r^{-2}\nabla_{\bbS^{n-1}} u\\
    &=\sum_{\ell\geq \ell_0}\sum_{m\in \calM_\ell}f_{\ell,m}r^{1-n-\ell}\left((2-n-\ell)Y_{\ell,m}(\omega)\hat{r}+r^{-1}\nabla_{\bbS^{n-1}} Y_{\ell,m}(\omega)\right).
\end{align*}
Note that the expression for the gradient differs from the usual one $\nabla u = \frac{\partial u}{\partial r}\hat{r}+\frac{1}{r}\nabla_{\bbS^{n-1}}u\hat{\theta}$ due to the normalization of $|\hat{\theta}|$ to $1$.

We then may estimate
\begin{align*}
    \|\sum_{\ell\geq \ell_0}\sum_{m\in \calM_\ell}f_{\ell,m}r^{1-n-\ell}\nabla_{\bbS^{n-1}} Y_{\ell,m}(\omega)\|_{L^2(B_R^+)}^2
    &\leq C\sum_{\ell\geq\ell_0}\sum_{m\in \calM_\ell}\ell^{-1}|f_{\ell,m}|^2\int_{\bbS^{n-1}}|\nabla_{\bbS^{n-1}}Y_{\ell,m}(\omega)|^2\:d\omega\\
    &\leq C \|f\|_{H^{1/2}(\partial B_R)}^2
\end{align*}
using the identity
\begin{align*}
    \langle \nabla_{\bbS^{n-1}}Y_{\ell,m},\nabla_{\bbS^{n-1}}Y_{\ell',m'}\rangle=-\langle \Delta_{\bbS^{n-1}}Y_{\ell,m},Y_{\ell',m'}\rangle
\end{align*}
for the orthogonality of the vector fields $\nabla_{\bbS^{d-1}} Y_{\ell,m}$ as well as to pass to the last line.
Bounding
\begin{align*}
    \|\sum_{\ell\geq 0}\sum_{m\in\calM_\ell}f_{\ell,m}r^{1-n-\ell}(2-n-\ell)Y_{\ell,m}(\omega)\hat{r}\|_{L^2(B_R^c)}\leq C\|f\|_{L^2(\partial B_R)},
\end{align*}
as before, we find that
\begin{align*}
    \|\nabla u\|_{L^2(B_R^+)}\leq C\|f\|_{H^{1/2}(\partial B_R)}.
\end{align*}
To proceed, we compute that for a function of the form $h(r)Y_{\ell,m}(\omega)$,
\begin{align*}
    \nabla^2 (hY)&=h''Y_{\ell,m} dr\otimes dr +\left(h'(r)-\frac{h(r)}{r}\right)(dr\otimes dY_{\ell,m}+dY_{\ell,m}\otimes dr)\\
    &\quad +rh'Y_{\ell,m}g_{\bbS^{n-1}}+h\nabla^2_{\bbS^{n-1}}Y_{\ell,m}.
\end{align*}
Therefore,
\begin{align*}
    \|\nabla^2 u\|_{L^2(B_R^+)}&\leq C(\|\sum_{\ell\geq \ell_0}\sum_{m\in \calM_\ell}\ell^2f_{\ell,m}r^{-n-\ell}dr\otimes dr\|_{L^2(B_R^+)}
    +\|\sum_{\ell\geq \ell_0}\sum_{m\in \calM_\ell}\ell f_{\ell,m}r^{1-n-\ell}dr\otimes dY_{\ell,m}\|_{L^2(B_R^+)}\\
    &\quad+\|\sum_{\ell\geq \ell_0}\sum_{m\in \calM_\ell}\ell f_{\ell,m}r^{2-n-\ell}Y_{\ell,m}g_{\bbS^{n-1}}\|_{L^2(B_R^+)}+\|\sum_{\ell\geq \ell_0}\sum_{m\in \calM_\ell} f_{\ell,m}r^{2-n-\ell}\nabla^2_{\bbS^{n-1}}Y_{\ell,m}\|_{L^2(B_R^+)}).
\end{align*}
The first three terms are estimated as before, whereas for the last term 
we write
\begin{align*}
    &\|\sum_{\ell\geq \ell_0}\sum_{m\in \calM_\ell} f_{\ell,m}r^{2-n-\ell}\nabla^2_{\bbS^{n-1}}Y_{\ell,m}\|_{L^2(B_R^+)}^2\leq\\
    &\quad\int\limits_0^\infty \left\|\nabla^2_{\bbS^{n-1}}\left(\sum_{\ell\geq \ell_0}\sum_{m\in \calM_\ell} f_{\ell,m}r^{2-n-\ell}Y_{\ell,m}\right)\right\|_{L^2(\bbS^{n-1})}^2 r^{n-1}\:dr
\end{align*}
and then apply elliptic regularity on $\bbS^{d-1}$:
\begin{align*}
&\int\limits_0^\infty \left\|\nabla^2_{\bbS^{d-1}}\left(\sum_{\ell\geq \ell_0}\sum_{m\in \calM_\ell} f_{\ell,m}r^{2-n-\ell}Y_{\ell,m}\right)\right\|_{L^2(\bbS^{n-1})}^2 r^{d-1}\:dr\leq\\
&\quad\int\limits_0^\infty \left\|\sum_{\ell\geq \ell_0}\sum_{m\in \calM_\ell} f_{\ell,m}r^{2-n-\ell}(1+\Delta_{\bbS^{n-1}})Y_{\ell,m}\right\|_{L^2(\bbS^{n-1})}^2r^{n-1}\:dr.
\end{align*}
Integrating in $r$ and using the orthogonality of the spherical harmonics, we conclude that
\begin{align*}
    \|\sum_{\ell\geq \ell_0}\sum_{m\in \calM_\ell} f_{\ell,m}r^{2-n-\ell}\nabla^2_{\bbS^{n-1}}Y_{\ell,m}\|_{L^2(B_R^+)}^2&\leq C\sum_{\ell\geq\ell_0}
    \sum_{\calM_\ell}\ell^{-1}|f_{\ell,m}|^2\|(1+\Delta_{\bbS^{n-1}} )Y_{\ell,m}\|_{L^2(\bbS^{n-1})}^2\\
&\leq C\|f\|_{H^{3/2}(\partial B_R)}^2.
\end{align*}
Having estimated $\|u\|_{L^2}$, $\|\nabla u\|_{L^2}$, and $\|\nabla^2 u\|_{L^2}$ in terms of $\|f\|_{H^{3/2}(\partial B_R)}$, the proof is now complete.
\end{proof}

\begin{lemma}
\label{lem:RnHelmholtz}
    Fix $E\in \bbC\setminus [0,\infty)$. Let $f\in H^{3/2}(\partial B_R)$ be given by
    \begin{align*}
        \sum_{\ell\geq 0}\sum_{m\in\calM_\ell}f_{\ell,m}Y_{\ell,m}. 
    \end{align*}
    Then the unique $u\in H^2(B_R^-)$ solving
    \begin{align}
    \label{eq:innerDeltaBVPHelmholtz}
    \begin{cases}
     (-\Delta-E) u=0\\
         T^- u=f,
    \end{cases}
    \end{align}
    is given by
    \begin{align*}
        u(r\omega)=\sum_{\ell\geq \ell_0}\sum_{m\in \calM_\ell}f_{\ell,m}\left( \frac{r}{R} \right)^{1-n / 2} \frac{K_{n/2+\ell-1}(z r)}{K_{n/2+\ell-1}(z R)}Y_{\ell,m}(\omega),
    \end{align*}
    where $z^2=-E$ with $\Re z>0$ and $K_\nu$ is the modified Bessel function of the second kind.
    Similarly, the unique $u\in H^2(B_R^+)$ solving
    \begin{align}
    \label{eq:outerDeltaBVPHelmholtz}
    \begin{cases}
     (-\Delta-E) u=0\\
         T^+ u=f,
    \end{cases}
    \end{align}
    is given by
    \begin{align*}
     u(r\omega)=\sum_{\ell\geq \ell_0}\sum_{m\in \calM_\ell}f_{\ell,m}\left( \frac{r}{R} \right)^{1-n / 2} \frac{I_{n/2+\ell-1}(z r)}{I_{n/2+\ell-1}(z R)}Y_{\ell,m}(\omega),
    \end{align*}
    where $I_\nu$ is the modified Bessel function of the second kind.
\end{lemma}
\begin{proof}
    As in the proof of Lemma~\ref{lem:exteriorZeroRn}, each $u_{\ell,m}$ satisfies the ODE 
   \begin{align*}
    r^2u''_{\ell,m}(r)+(n-1)ru'(r)-\ell(\ell+n-2)u(r)=0.
    \end{align*}
    Using the change of variables
    \begin{align*}
            h(r)=r^{\frac{n}{2}-1}u_{\ell,m}(r),
    \end{align*}
    it is easy to see that basis of solutions is given by $\{r^{1-n/2}I_{n/2+\ell-1}(z r),r^{1-n/2}K_{n/2+\ell-1}(zr)\}$
        From \cite[10.17.3-4]{AS}, we have that as $r\to \infty$, for some constants $C_{1,2}>0$
        \begin{align*}
           &I_\nu (r)=C_1\frac{e^{r}}{(2\pi r)^{\frac{1}{2}}}+O(\frac{1}{r^{\frac{1}{2}}})\\
           &K_\nu ( r)=C_2\frac{\sqrt{\pi}e^{-r}}{(2 r)^{\frac{1}{2}}}+O(\frac{1}{r^{\frac{1}{2}}})
        \end{align*}
    and therefore the solution $B_R^+$ must be proportional to $r^{1-n/2}K_{n/2+\ell-1}(zr)$. On the other hand, from \cite[10.30(i)]{AS}, we have that as $r\to 0$, for some constants $C_{1,2}>0$
        \begin{align*}
           &I_\nu (r)=C_1(\frac{r}{2})^{\nu}+O(r^{\nu+1})\\
           &K_\nu ( r)=C_2(\frac{r}{2})^{-\nu}+O(r^{-\nu+1}),
        \end{align*}
    and therefore the only solution which is $L^2$ near $0$ must be proportional to $r^{1-n/2}I_{n/2+\ell-1}(zr)$. Plugging in the boundary condition, the result now follows.
\end{proof}

\bibliographystyle{amsplain}
\bibliography{bib}

\end{document}